\documentclass{amsart}

\usepackage{amsthm}
\usepackage{amssymb}
\usepackage{amsmath}
\usepackage{amscd}
\usepackage{epic}
\usepackage{eepic}

\newcommand{\labell}[1] {\label{#1}}
\newcommand    {\comment}[1] {}

\newtheorem {Theorem}   {Theorem}

\numberwithin{Theorem}{section}

\newtheorem {Lemma}[Theorem]    {Lemma}         
\newtheorem {Proposition}[Theorem]{Proposition}  
\theoremstyle{definition}
\newtheorem{Definition}[Theorem]{Definition}
\theoremstyle{remark}
\newtheorem{Remark}[Theorem]{Remark}
\newtheorem{Example}[Theorem]{Example}

\newtheorem {Corollary}[Theorem]{Corollary}  
\newtheorem {Question}[Theorem]{Question}

\def	\HA	{\operatorname{HA}}

\def    \R      {{\mathbb R}}
\def    \reals  {{\mathbb R}}

\def    \C      {{\mathbb C}}

\def	\Cx	{{{\mathbb C}^\times}}
\def    \CP     {{\mathbb C}{\mathbb P}}
\def    \Z      {{\mathbb Z}}
\def    \T      {{\mathbb T}}
\def    \N      {{\mathbb N}}
\def    \I      {{\mathbb I}}

\def    \AA     {{\mathcal A}}
\def	\calC	{{\mathcal C}}
\def    \EE     {{\mathcal E}}

\def    \LL     {{\mathcal L}}
\def    \MM     {{\mathcal M}}

\def    \MMs    {{\mathcal M}^{{\mathit smooth}}}
\def    \OO     {{\mathcal O}}
\def    \PP     {{\mathcal P}}

\def	\Nu	{{\mathcal V}}

\def    \inv    {^{-1}}

\def    \im     {{\operatorname{im}}\,}

\def    \ol     {\overline}

\def    \g      {{\mathfrak g}}

\def    \h      {{\mathfrak h}}

\def    \k      {{\mathfrak k}}

\def    \ssminus        {{\smallsetminus}}

\def    \p      {\partial}

\def    \id     {\operatorname{id}}

\def    \codim   {{\operatorname{codim}}}

\def    \inv    {^{-1}}
\def    \ssminus        {\smallsetminus}

\def    \ol     {\overline}

\def    \tPsi   {\tilde{\Psi}}

\def	\zbar	{{\bar{z}}}
\def \liminv {\lim\limits_{\longleftarrow}}
\def \liminvk {\sideset{}{^{(k)}}\lim\limits_{\longleftarrow}}

\begin{document}


\setlength{\smallskipamount}{6pt}
\setlength{\medskipamount}{10pt}
\setlength{\bigskipamount}{16pt}





\title[Assignments and Abstract Moment Maps]
{Assignments and Abstract Moment Maps}

\author[Viktor Ginzburg]{Viktor L. Ginzburg}
\author[Victor Guillemin]{Victor Guillemin}
\author[Yael Karshon]{Yael Karshon}

\address{V. Ginzburg: Department of Mathematics, U. C. Santa Cruz,
Santa Cruz, CA 95064, USA}
\email{ginzburg@math.ucsc.edu}

\address{V. Guillemin: Department of Mathematics, MIT, Cambridge, 
MA 02139-4307, USA}
\email{vwg@math.mit.edu}

\address{Y. Karshon: The Hebrew University of Jerusalem, Institute of
Mathematics, Giv'at Ram, 91904 Jerusalem, Israel}
\email{karshon@math.huji.ac.il}

\date{\today}

\thanks{The work of all three authors was partially supported by the NSF
and by the Israel--US Binational Science Foundation.}

\begin{abstract}
Abstract moment maps arise as a generalization of genuine moment
maps on symplectic manifolds when the symplectic structure is
discarded, but the relation between the mapping and the action is kept.
Particular examples of abstract moment maps had been used in Hamiltonian
mechanics for some time, but the abstract notion originated in the study of 
cobordisms of Hamiltonian group actions.

In this paper  we answer the question of existence of a (proper) abstract 
moment map for a torus action and give a necessary and
sufficient condition for an abstract moment map to be associated
with a pre-symplectic form. This is done by using the notion of an
assignment, which is a combinatorial counterpart of an abstract moment map. 

Finally, we show that the space of assignments fits as the zeroth 
cohomology in a series of certain cohomology spaces associated with a torus 
action on a manifold. We study the resulting ``assignment cohomology" theory.
\end{abstract}

\maketitle

\bigskip

\section{Introduction}
\labell{sec:intro}

Abstract moment maps arise as a generalization of genuine moment maps 
on symplectic manifolds. The essence of their definition is that
the symplectic structure is discarded but the relation between the 
mapping and the action is kept intact. To be precise, an abstract moment 
map on a $G$-manifold $M$
is an equivariant mapping $\Psi\colon M\to\g^*$ satisfying the following
constancy condition: for every Lie algebra element $\xi\in\g$, the component
$\left<\Psi,\xi\right>$ is locally constant on the set $\{\xi_M=0\}$  
where the action generating vector field $\xi_M$ vanishes.

Moment maps on symplectic manifolds are among examples
of abstract moment maps. In general, however, an abstract moment map is 
not associated with a symplectic form or even a closed two-form. An 
abstract moment map is an additional structure on a $G$-manifold, 
and a given $G$-manifold admits many abstract moment maps. 
Thus $G$-manifolds
equipped with abstract moment maps occupy an intermediate place between 
pure $G$-manifolds and symplectic manifolds with Hamiltonian $G$-actions.

The goal of this paper is to find a relationship between
$G$-manifolds, $G$-manifolds equipped with abstract moment maps,
and Hamiltonian $G$-spaces. This is done by using a new notion, the notion
of an \emph{assignment} comprising certain combinatorial data
extracted from an abstract moment map. 
An assignment
should be thought of as a combinatorial counterpart of an abstract 
moment map. For a torus action, an assignment is a function
from the set of orbit type strata to the dual spaces to
the Lie algebras of the stabilizers. (Thus a manifold with a finite 
orbit type stratification has only a finite--dimensional space of 
assignments.)

First we address the existence and uniqueness question for abstract moment
maps with a fixed assignment. We prove that, for a torus action, every 
assignment is associated with an abstract moment map. Furthermore,
two abstract moment maps give rise to
the same assignment if and only if they differ, roughly speaking,
by a Hamiltonian moment map arising from an exact two-form.

For some problems concerning non-compact manifolds, it is important to 
consider abstract moment maps which are proper. 
(Non-compact $G$-manifolds with proper abstract moment maps share some
of the appealing properties of compact $G$-manifolds. 
See \cite{Ka,GGK:lerman,GGK:book}.) 
We show that for a given assignment, an abstract moment map can be chosen 
proper if the assignment is proper or, to be more precise, ``polarized"; 
see Section \ref{subsec:polar}.

Then we use the notion of an assignment to answer the question 
whether a given abstract map is associated with a two-form.

Moment maps associated with true symplectic forms 
must additionally satisfy some non-degeneracy 
requirements. These requirements are analyzed in \cite{GGK:book}. 
On the other hand, moment maps on
Poisson manifolds (see, e.g.,  \cite{wein-lect})
are not in general abstract moment maps because they
do not have to be locally constant on the fixed point set.

The paper is organized as follows. In Section \ref{sec:def-examples}
we recall the definition of abstract moment maps and illustrate it
by a number of examples.
In Section \ref{sec:existence} we give a necessary and
sufficient condition (in terms of assignments) for a $G$-manifold to 
admit a (polarized) abstract moment map. 
In Section \ref{sec:exact} we show that two abstract moment maps
with the same assignment differ by  one
which is associated with
an exact two-form (or, to be more precise, with a one-form). 
We call such abstract moment maps \emph{exact}.
Section \ref{sec:forms} is devoted to the question of which abstract moment
maps are Hamiltonian, i.e., associated with closed two-forms. 
We show that every abstract moment map
is locally Hamiltonian. Globally, there is an obstruction, which is stated
again in terms of assignments.
In Section \ref{sec:proof} we prove a technical theorem on which
the results of Sections \ref{sec:exact} and \ref{sec:forms} heavily rely:
an abstract moment map on a linear representation
is exact if and only if it vanishes at the origin.
Finally, we show that the space of assignments and some of its 
generalizations fit as the zeroth cohomology in a series of certain 
cohomology spaces associated with a $G$-manifold. This cohomology is
introduced and studied in Section \ref{sec:cohomology}.

Abstract moment maps were introduced in \cite{Ka} to study geometric
equivariant $G$-cobordisms and to state and prove
the cobordism linearization theorem.
This theorem, whose earliest version was given in \cite{GGK1} 
(see also \cite{GZ} for important related work), 
asserts that under certain hypotheses a manifold with a torus action
is equivariantly cobordant to the disjoint union of the linear isotropy 
representations at the fixed points.\footnote{
A number of applications of the linearization theorem  
are outlined in \cite{GGK1} and \cite{Ka}. 
See \cite{MW1} and \cite{MW2} for some new applications, 
and see \cite{GGK:book} and \cite{GGK:lerman} for further developments.
}

One of the main conceptual difficulties in the formulation of this 
theorem is to find a notion of non-compact cobordism which
would not render every compact manifold cobordant to the empty set.
(This is the central problem arising when non-compact manifolds are introduced
in a cobordism theory: all compact manifolds may become
cobordant to each other and so to zero.)
The notion of an abstract moment map provides a solution 
to this problem: two $G$-manifolds
equipped with proper abstract moment maps are cobordant if there exists
a $G$-cobordism between them and a proper abstract moment map on it extending
those on the boundary. This definition leads to a non-trivial cobordism
theory in which the theory of compact (geometric) $G$-cobordisms is 
embedded.\footnote{
Strictly speaking, this is true for geometric stable complex $G$-cobordisms. 
See \cite{GGK:lerman} and references therein.
}
The non-compact theory appears to be in some sense simpler than the compact one.
The reason is that the non-compact theory has a 
well-understood set of generators and probably fewer relations than the compact one.
(See \cite{Ka} and \cite{GGK:book}.)

We feel, however, that, as some examples in Sections 
\ref{sec:def-examples} and \ref{sec:cohomology} indicate, 
abstract moment maps and assignments may have uses beyond
those connected with geometric equivariant $G$-cobordisms.

\subsection*{Notation and conventions} 
Throughout this paper, $M$ is a $G$-manifold ($C^{\infty}$-smooth),
with $G$ being a torus, except in some rare cases where $G$ is allowed 
to be a more general compact Lie group. 
As usual, $\g$ denotes the Lie algebra of $G$ and $\xi_M$
is the vector field induced by the action of $\xi\in\g$ on $M$. The
stabilizer of $x\in M$ is denoted $G_x$ and the fixed point set of $G$ on
$M$ by $M^G$.
All ordinary and equivariant cohomology 
groups are assumed to have real coefficients unless specified otherwise.

We consider only abstract moment maps for
torus actions. Many (but not all) of our results should extend
to proper actions of other Lie groups. 

\subsection*{Acknowledgment}
The authors wish to thank Phil Bradley, Emmanuel Farjoun, Debra Lewis,
Assaf Libman, David Metzler, Avishay Vaaknin, and Yuli Rudyak for fruitful 
discussions and comments.

\section{Abstract moment maps}
\labell{sec:def-examples}

Let us first recall some standard facts about smooth group actions.
Let a Lie group $G$ act smoothly on a manifold $M$.
Assume that $G$ is compact or, more generally, that the action is proper
(that is, that the map $(a,m) \mapsto (a \cdot m , m)$ from $G \times M$
to $M \times M$ is proper). 
Each Lie algebra element $\xi \in \g$ gives rise to a vector field $\xi_M$
on $M$. The \emph{stabilizer} of a point $m \in M$ is the group
$\{a \in G \ | \ a \cdot m = m \}$.  The Lie algebra of this group
is equal to $\{ \xi \in \g \ | \ \xi_M|_m = 0 \}$
and is called the \emph{infinitesimal stabilizer} of $m$.
For each subgroup $H \subseteq G$, the connected components of the set 
of points whose stabilizer is conjugate to $H$ (hence equal to $H$,
if $G$ is Abelian) are smooth sub-manifolds of $M$. These connected
components are the
\emph{orbit type strata} of $M$. They are partially ordered:
$X \preceq Y$ if and only if the stratum $X$ is contained in the closure
of the stratum $Y$. Similarly, the connected components of the sets 
of points whose infinite stabilizers are conjugate to given sub-algebras
$\h \subseteq \g$ (hence equal to $\h$, if $G$ is Abelian)
form the \emph{infinitesimal obit type strata} of $M$. 
These are, too, partially ordered by inclusions of closures of strata.
The infinitesimal orbit type stratification is more coarse than the
orbit type stratification, because two points $x$ and $y$ can have
different stabilizers but the same infinitesimal stabilizer.
In what follows we will mainly work with
the infinitesimal orbit type stratification because this stratification
is more suitable for the goals of this paper than the orbit type 
stratification. 

Next, let us recall the definition of abstract moment maps, 
as given in \cite{Ka}. 
For a map $\Psi \colon M \to \g^*$, we denote by $\Psi^\h$ 
or $\Psi^H$ the composition of 
$\Psi$ with the natural projection $\g^* \to \h^*$, where $\h$ is the
Lie algebra of $H$.  Similarly, for any Lie algebra element $\xi \in \g$, 
we denote by $\Psi^\xi \colon M \to \R$ the $\xi$th component of $\Psi$, 
i.e., $\Psi^\xi=\left<\Psi,\xi\right>$.

\begin{Definition} \labell{def:moment}
An \emph{abstract moment map} on $M$ is a smooth map 
$ \Psi \colon M \to \g^* $ with the following properties:
\begin{enumerate}
\item
$\Psi$ is $G$-equivariant, and
\item
for any subgroup $H$ of $G$, the map
$\Psi^H \colon M \to \h^* $ is locally constant on the submanifold $M^H$
of points fixed by $H$.
\end{enumerate}
\end{Definition}

\begin{Remark}
For the second requirement to hold, it is enough to assume that 
for any Lie algebra element $\xi \in \g$, 
the function $\Psi^\xi$ is locally constant on the set of zeros
of the corresponding vector field $\xi_M$.
If $G$ is compact, it is enough to demand the requirement
for circle subgroups of $G$.
\end{Remark}

In this paper we mainly consider the case where $G$ is a torus and we 
often focus on abstract moment maps that are proper.

\begin{Example} \labell{ex:zero}
The constant function zero is an abstract moment map. 
It is proper if and only if $M$ is compact.
\end{Example}

\begin{Example}
If the fixed point set $M^G$ has a non-compact component, 
$M$ does not admit a proper abstract moment map.
\end{Example}

\begin{Example}
Let $G$ be the circle group, and let us identify $\g^*$
with $\R$. Then an abstract moment map is a real valued 
invariant function that is constant on each connected component 
of the fixed point set.  In particular, if the set of fixed points is 
discrete, any invariant function is an abstract moment map. 
\end{Example}

\begin{Example} \labell{exam:presymplectic}
Recall that a Hamiltonian $G$-space is a triple $(M,\omega,\Psi)$,
where $M$ is a $G$-manifold, $\omega$ is a closed invariant two-form 
(which in some contexts -- not here -- is required to be symplectic),
and $\Psi$ is a moment map, i.e., a $G$-equivariant function 
$\Psi \colon M \to \g^*$, such that Hamilton's equation,
\begin{equation} \labell{Hamilton}
 \iota(\xi_M) \omega = - d\Psi^\xi \quad \text{for all $\xi \in \g$},
\end{equation}
holds. Then $\psi$ is an abstract moment map.
\end{Example}

If equation \eqref{Hamilton} holds, we say that $\omega$ 
is compatible with $\Psi$, or that $\Psi$ is associated with $\omega$. 
An abstract moment map associated with some two-form will be called 
a \emph{Hamiltonian} moment map.

\begin{Example} \labell{ex:one-form}
Let a Lie group $G$ act on a manifold $M$
and let $\mu$ be any invariant one-form.
Then the function $\Psi \colon M \to \g^*$ defined by
\begin{equation} \labell{moment one form} 
\Psi^\xi = \mu(\xi_M)
\end{equation}
is an abstract moment map.
Moreover, for each $H \subset G$, the function $\Psi^H$ vanishes on $M^H$.
\end{Example}

\comment{\ In the book, expand the following paragraph:}

An abstract moment map that arises by Equation \eqref{moment one form}  
is called \emph{exact}. 
A compatible two-form is then given by $\omega = d\mu$. 
Many of ``classical'' moment
maps, e.g., the canonical moment map on the cotangent bundle, are exact.
Also, in the pre-quantization of a Hamiltonian action, 
the pullback of the moment map to the pre-quantum circle bundle 
is an exact abstract moment map.

The advantage of working with exact moment maps over Hamiltonian ones
is  that if $\Psi_0$ and $\Psi_1$
are exact moment maps then so is $(1 - \rho) \Psi_0 + \rho \Psi_1$ 
for any smooth function $\rho$.

\begin{Remark}
Recall that for any Lie group $G$ an equivariant differential two-form on a $G$-manifold $M$
is a formal sum 
\begin{equation} \labell{omegaPsi}
\omega + \Psi,
\end{equation}
where $\omega$ is an invariant two-form on $M$ and $\Psi$ is a smooth 
equivariant function from $M$ to $\g^*$.  
The equivariant form \eqref{omegaPsi} 
is said to be equivariantly closed if and only if 
it satisfies \eqref{Hamilton};
it is said to be equivariantly exact if and only if there exists 
an invariant one-form $\mu$ such that $\omega = d\mu$ and 
$\Psi^\xi = \mu(\xi_M)$ for all $\xi \in \g$.
The second equivariant cohomology, denoted $H^2_G(M)$,
is the quotient of the space of equivariantly closed equivariant two-forms
by the subspace of those that are equivariantly exact.
\end{Remark}

\begin{Example} \labell{pull-back}
The pull-back of an abstract moment map is an abstract moment map. More
precisely, let $f\colon N\to M$ be an equivariant map of $G$-manifolds
and let $\Psi$ be an abstract moment map on $M$. Then $f^*\Psi=\Psi\circ f$
is an abstract moment map on $N$; the map $f^*\Psi$ is proper, provided that
$f$ and $\Psi$ are proper.

For instance, following \cite{SLM} and \cite{Le}, 
consider a $G$-manifold $Q$ 
and denote by $J\colon T^*Q\to \g^*$ the canonical moment map:
$J^\xi(\mu)=\mu(\xi_Q(x))$, where $\mu\in T^*_xQ$.
The action map $F\colon Q\times \g\to TQ$ is defined as $F(x,\xi)=\xi_Q(x)$.
Consider a Lagrangian on $Q$ with Legendre transformation  
$\LL \colon TQ\to T^*Q$. Then 
$$
\I=(\LL F)^*J\colon Q\times \g\stackrel{F}{\to} TQ\stackrel{\LL}{\to} T^*Q
\stackrel{J}{\to} \g^*
$$
is an exact abstract moment map. For example, assume that $\LL$ arises from a 
Riemannian metric 
$\left<~,~\right>$ on $Q$ so that $\LL(v)= \left<v,\cdot\right>$ for
a tangent vector $v$. Then
$\I^\zeta(x,\xi)=\left<\xi_Q(x),\zeta_Q(x)\right>$.
The map $\I$, called the locked momentum map, is used in the analysis of 
relative equilibria. (See \cite{SLM} and \cite{Le}.) Note that in general 
$Q\times \g$ is not a 
symplectic manifold and $G$, in this example, does not have to be commutative.
\end{Example}

\begin{Remark}
In view of Example \ref{exam:presymplectic}, it is worth
pointing out that a moment map on a Poisson manifold 
(see, e.g., \cite{wein-lect})
may not be an abstract moment map even when the Poisson structure
is preserved by the action. The reason is that, since a moment
map is defined only up to addition of Casimir functions, and since on a 
Poisson manifold Casimir functions often exist in abundance, 
a moment map on a Poisson manifold may not be constant on the fixed 
point set. 
\end{Remark}

\section{Existence of abstract moment maps}
\labell{sec:existence}

Every manifold with a $G$-action admits an abstract moment map:
the zero map. This map is never proper unless the manifold 
is compact. In this section we answer the question of when a $G$-manifold 
admits a proper (in fact, polarized, see Definition \ref{polarized}) 
abstract moment map.

A necessary condition for an action to admit a proper abstract moment map 
$\Psi$ is that each component of the fixed point set be compact. 
(Recall that $\Psi$ is constant on each such component.) 
Is this condition sufficient?
Moreover, does a (proper) abstract moment map exist with prescribed values
at the fixed points?

\subsection{Existence of abstract moment maps for circle actions}

Answers to the above questions take a particularly simple and attractive
form when $G$ is a circle, when abstract moment maps are simply
$G$-invariant functions that are constant on the connected components
of the fixed point set.

\begin{Theorem} \labell{thm:circle-exist}
Let $G$ be a circle acting on $M$, and let $\psi \colon M^G\to \reals$ be
a locally constant function.
\begin{enumerate}
\item There exists an abstract moment map $\Psi\colon M\to\reals$ with
$\Psi|_{M^G}=\psi$.
\item Assume that $\psi$ is proper and bounded from below. Then $\Psi$
can be chosen to be proper and bounded from below.
\end{enumerate}
\end{Theorem}

\begin{Remark}
In other words, if $G$ is the circle group, we can prescribe the values
of an abstract moment map on the connected components of $M^G$ completely 
arbitrarily. If $M^G$ is compact, the condition of the second assertion
is satisfied automatically, and every locally constant function 
on $M^G$ extends to a proper abstract moment map.
\end{Remark}

\begin{proof}
The theorem follows from the following two facts, 
applied to $X = M^G$ and $f = \psi$.

\begin{enumerate}
\item \labell{part 1}
Let $X \subset M$ be a closed sub-manifold and $f \colon X \to \R$ a smooth
function.
Then there exists a smooth function $F \colon M \to \R$ such that $F|_X=f$. 
Moreover, if $f$ is bounded from below, $F$ can be chosen to be
bounded from below too, and if $f$ is proper and bounded from below,
$F$ can be chosen to be proper and bounded from below.

\item \labell{part 2}
Let $F \colon M \to \R$ be proper and bounded from below.
Then the average $\ol{F}$ of $F$ by a compact group action is also 
proper and bounded from below.
\end{enumerate}

Let us prove the first fact. Fix a tubular neighborhood $U$ of $X$ 
in $M$, and let
$\pi \colon U \to X$ be a smooth projection which extends to a proper
map from the closure $\ol{U}$ to $X$.
Let $\rho$, $1-\rho$ be a smooth partition of unity subordinate
to the covering of $M$ by the two open sets $U$ and $M \ssminus X$.
Pick a smooth function $\varphi \colon M \to \R$ which is proper and 
bounded from below (see, e..g., \cite{GP}, Chapter 1, Section 8).
Then $F = \rho f + (1 - \rho) \varphi$ has the desired properties.

To prove the second fact, notice that ${\ol{F}}^{-1}([-a,a])$
is contained in $G\cdot F^{-1}([-a,a])$, which is the image of
the compact set $G\times F^{-1}([-a,a])$ under the continuous action
mapping $G\times M\to M$. 
\end{proof}

\begin{Remark}
In Theorem \ref{thm:circle-exist} it is not true that if $\psi$ 
is just proper (but not bounded) then $\Psi$ can be chosen to be proper. 
In general, a proper map on a closed submanifold $X$ of $M$ might not
extend to a proper map on $M$.  For instance, the function $f(0,y) = y$
on the $y$-axis does not extend to a continuous proper function $F$ from
$\R^2$ to $\R$. A similar counterexample involving abstract moment maps
is given below.
\end{Remark}

\begin{Example}
Let $M$ be obtained by the following plumbing construction:
$$ M = \Z \times S^2 \times D^2 / \sim,$$
where $S^2 = \{ (x,y,z) \in \R^3 \ | \ x^2 + y^2 + z^2 =1 \}$,
$D^2 = \{ (u,v) \in \R^2 \ | \ u^2 + v^2 < \epsilon^2 \}$,
and 
$(n,x,y,\sqrt{1-x^2-y^2},u,v) \sim (n+1,u,v,-\sqrt{1-u^2-v^2},x,y)$
for all $n$. Take the diagonal circle action:
$$ e^{i\theta} \cdot [n,x,y,z,u,v] = [n,x',y',z,u',v']$$
where 
$$ \left[ \begin{array}{cc} x' & u' \\ y' & v' \end{array} \right]
 = \left[ \begin{array}{cc} 
\cos\theta & \sin\theta \\ -\sin\theta & \cos\theta \end{array} \right]
  \left[ \begin{array}{cc} x & u \\ y & v \end{array} \right].
$$
The function
\begin{equation} \labell{prop not pol}
 \psi([n,0,0,1,0,0]) = n
\end{equation}
is a locally constant function on the fixed point set and is proper,
but it does not extend to a proper function $\Psi\colon M \to \R$.
(The function $\psi$ extends to a (non-proper) Hamiltonian moment map,
 for a closed two-form whose pullback to each 
$\{ n \} \times S^2 \times \{ 0 \}$ is non-negative and has total area one.)
\end{Example}

\subsection{Assignments}
Let us now investigate more closely the question of existence
of an abstract moment map for an action of a torus whose dimension
is greater than one.  Theorem \ref{thm:circle-exist}
is no longer true in this case:

\begin{Example}
\labell{exam:S2xS2}
Let $\Psi = (\Psi_1,\Psi_2)$ be an abstract moment map on 
$M=S^2 \times S^2$ with $G = S^1 \times S^1$ acting
by rotating each of the two factors.
Then $\Psi$ must send the four fixed points
to the corners of a rectangle in $\R^2$ whose sides are parallel to the
axes. Thus the values $\Psi(M^G)$ cannot be assigned arbitrarily.
\end{Example}

Moreover, the abstract moment maps might not even separate the components
of the fixed point set:

\begin{Example} \labell{S4}
Let  $S^4$  be the unit sphere in $\C\times\C\times\R$,  
and let  $G=S^1\times S^1$ act on it by rotating 
each of the the first two factors.  
There are two fixed points: the North Pole and the South Pole.  
The set of points fixed by the first  $S^1$  is connected 
and contains both poles. The same is true for the set of points
fixed by the second  $S^1$.  Consequently, any abstract moment map on  $S^4$
must have the same value at the poles.
\end{Example} 

These examples stress the role of the orbit type strata other than the
components of the fixed point set.

For a stratum $X$ we denote by $\g_X$ the infinitesimal stabilizer of
any of its points.

Suppose that $\Psi \colon M \to \g^*$ is an abstract moment map. 
Then for each infinitesimal orbit type stratum $X$ in $M$, the map
$\Psi$ followed by the projection $\g^* \to \g_X^*$ gives an element 
$A(X)$ of $\g_X^*$. 
An important observation is that the existence question for an abstract 
moment map is equivalent to the existence question for such an assignment,
$X \mapsto A(X)$.
We make this precise in Theorems \ref{thm:assignment} and 
\ref{thm:assignment proper}, which rely on the following definition and 
example.

\begin{Definition}
\labell{def-assign}
An \emph{assignment} is a function $A$ that associates to each infinitesimal
orbit type stratum $X$ in $M$ an element $A(X)$ of $\g_X^*$
and that satisfies the following compatibility condition:
if $X$ is contained in the closure of $Y$
then $A(Y)$ is the image of $A(X)$ under the restriction map
$\g_X^* \to \g_Y^*$.
The linear space of all assignments on $M$ is denoted by $\AA(M)$.
\end{Definition}

In Section \ref{sec:cohomology} we discuss assignments in a broader,
more abstract, context.

\begin{Example}
\labell{exam:moment-assign}
Let $\Psi \colon M \to \g^*$ be an abstract moment map.
Then $A(X) = \Psi^{\g_X}(X)$ is an assignment. The assignment $A$
and the moment map $\Psi$ are said to be associated with each other. 
If the abstract moment map is exact, i.e., $\Psi$ arises from a 
one-form $\mu$ so that $\Psi^\xi = \mu(\xi_M)$, the corresponding 
assignment is zero.
\end{Example}

\begin{Example}
\labell{exam:assign-circle}
When $G$ is the circle group, an assignment simply associates a real number
to each component of the fixed point set.
Thus, in this case, $\AA(M)=(\g^*)^{\pi_0(M^G)}$.
\end{Example} 

\begin{Example} \labell{CP2-T2}
Consider the action of the two-dimensional torus $G=S^1 \times S^1$ 
on $M=\CP^2$ given by $(t_1,t_2)[z_0:z_1:z_2]=[z_0:t_1z_1:t_2z_2]$.
This action has three fixed points,
and every assignment $A$ is uniquely determined by its value at the
fixed points. There are, however, three relations between the values
$A|_{M^G}\in(\g^*)^3$, coming from the strata with one-dimensional
stabilizers. As a result, $\AA(M)$ is three-dimensional.
Geometrically, the assignment values at the fixed points are the vertices
of a triangle in $\R^2$ with two equal sides that are parallel 
to the coordinate axes.
\end{Example}

\begin{Example}
\labell{exam:product}
Let $M_1$ be a $G_1$-manifold and $M_2$ be a $G_2$-manifold. Then
$$
\AA(M_1\times M_2)=\AA(M_1)\oplus \AA(M_2)
$$
for the $G_1\times G_2$-action on $M_1\times M_2$.
\end{Example} 

\begin{Remark}
Replacing in Definition \ref{def-assign} the infinitesimal orbit type 
stratification by the orbit type stratification notion leads
to the same class of assignments $\AA(M)$. 
Namely, a function that associates to each orbit type stratum $X$
an element of $\g_X^*$ and that satisfies the compatibility condition
of Definition \ref{def-assign} is in fact constant on each infinitesimal
orbit type stratum, and, hence, is an assignment.
\end{Remark}

\begin{Example} \labell{toric:1}
Let $M$ be an $n$-dimensional complex manifold and let $G$ be 
an $n$-dimensional 
torus that acts on $M$. Suppose that each point of $M$ with stabilizer
$H \subseteq G$ has a neighborhood which is biholomorphic to a neighborhood
of the origin in $\C^n$ with an action of $H$ of the following form.
The $H$-action is obtained as the composition of an isomorphism 
$H\to (S^1)^{\dim H}$ with the $(S^1)^{\dim H}$-action on 
$\C^{\dim H} \times \C^{n - \dim H}$ which is standard on the first factor
and trivial on the second. (For instance, this is the case if $M$ is a 
toric manifold; see \cite{fulton} or \cite{audin}.)
It follows that for any stratum $X$, the natural map
$$ \g_X^* \stackrel{\oplus_Y \pi^X_Y}{\longrightarrow} 
   \bigoplus\limits_{\{ Y \ | \ X \preceq Y, \ \dim \g_Y=1 \}} \g_Y^* $$
is a linear isomorphism. Therefore, a moment assignment is determined
by its values on the strata $Y$ with $\dim \g_Y = 1$, and these values
can be prescribed arbitrarily. So for such $M$,
$$ 
\AA(M) =\bigoplus_Y\g_Y^*\cong \R^{\# \{Y \ | \ \dim \g_Y =1 \} }.
$$
\end{Example}

\begin{Example} \labell{toric:2}
Let $M$ be a K\"ahler toric manifold (also see Example \ref{toric:1})
with moment map $\Psi \colon M \to \g^*$. The image $\Psi(M)$ is a convex
polytope \cite{At:convexity,GS:convexity}. 
A convex polytope is stratified, with the strata being its open faces 
of various dimensions. 
The orbit type strata in $M$ are exactly the preimages in $M$ of the 
open faces in $\Phi(M)$, \cite{Del}.
As a consequence, the poset of (infinitesimal) orbit type strata of $M$
is isomorphic to the poset of faces of $\Psi(M)$.
Moreover, the stabilizers can be read from the faces
and vice versa: the affine plane spanned by the face $\Psi(X)$ 
is a shift of the annihilator in $\g^*$ of the Lie algebra $\g_X$.
The shifts are exactly given by the assignment 
$A(X) = \Psi^{\g_X} (X)$:
$$ \text{affine span} (\Psi(X)) 
   = \text{ preimage of $A(X)$ under $\g^* \to \g_X^*$ }.
$$
The polytope $\Psi(M)$, and hence the moment assignment $A$,
determine the manifold, the $G$-action, and the symplectic form
up to an equivariant symplectomorphism, \cite{Del},
and the equivariant K\"ahler structure on the strata, \cite{Gu}.
\end{Example}

\begin{Remark}
In Example \ref{toric:2} we saw that the moment assignment of a symplectic
toric manifold $M$ determines its moment polytope $\Psi(M)$.
Similarly, for a toric manifold with a closed invariant two-form
which may have degeneracies, the moment assignment determines
its \emph{twisted polytope} in the sense of \cite{KT:toric}.
\end{Remark}

An obvious, but important, fact is
\begin{Lemma} \labell{fix A}
Let $\Psi_0$ and $\Psi_1$ be abstract moment maps which have the same 
assignment, $A$. Then $(1-\rho) \Psi_0 + \rho \Psi_1$ is also an abstract
moment map with assignment $A$, for any invariant smooth function $\rho$.
\end{Lemma}

\begin{proof}
For any $H \subseteq G$, on every component $X$ of $M^H$, we have
$$
(1-\rho) \Psi_0^H + \rho \Psi_1^H = (1-\rho) A(X) + \rho A(X) \equiv A(X).
$$
\end{proof}

\begin{Remark}
\labell{rmk:non-abelian}
The definition of an assignment can be extended to non-commutative groups
$G$.   An assignment can then be defined as a function $x\mapsto
A(x)\in \g_x^*$ on $M$ such that the following conditions hold: 
\begin{itemize}

\item $A(g\cdot x)=Ad^*_g A(x)$ for all $x\in M$ and $g\in G$.

\item $A^\h$ is locally constant on the set $M^\h$ of points $x$ with 
$\h\subseteq \g_x$.
\end{itemize}
In the non-commutative case, as in Example \ref{exam:moment-assign}, 
an abstract moment map gives rise to an assignment.
\end{Remark}

\subsection{Existence of abstract moment maps for torus actions}
\labell{subsec:polar}
The relation between abstract moment maps and assignments is expressed in
the following

\begin{Theorem} \labell{thm:assignment}
Let $M$ be a manifold with a $G$ action.
Let $A: X \mapsto A(X)$ be an assignment. Then there exists an abstract 
moment map $\Psi \colon M \to \g^*$ which is associated with $A$,
i.e., such that $\Psi^{\g_X}(X) = A(X)$ in $\g_X^*$ 
for every orbit type stratum $X$.
\end{Theorem}

\begin{proof}
Let $m$ be a point in $M$, let $\h$ be the infinitesimal stabilizer of $m$,
and let $A(m) \in \h^*$ be the element assigned to the orbit type
stratum containing $m$.
Let $\Psi_m \in \g^*$ be any element whose projection to $\h^*$ is $A(m)$.
Pick an open neighborhood $U_m$ of the orbit $G \cdot m$
which equivariantly retracts to the orbit. The constant function
$\Psi_m$ is an abstract moment map on $U_m$ whose assignment is $A|_{U_m}$.
Choose an invariant partition of unity $\{\rho_j\}$ subordinate to the
covering of $M$ by the open subsets $U_m$, with the support of $\rho_j$
contained in the open set $U_{m_j}$. The convex combination
$\Psi = \sum \rho_j \Psi_{m_j}$ is an abstract moment map;
this follows from Lemma \ref{fix A}, applied to open subsets of the
manifold.
\end{proof}

On a non-compact manifold, it is sometimes required that 
an abstract moment map be proper. (See \cite{Ka,GGK:lerman,GGK:book}.)
In fact, we often need a component of $\Psi$ to be proper and bounded 
from below.

\begin{Definition} \labell{polarized}
Let $\eta \in \g$ be a Lie algebra element.
A function $\Psi \colon M \to \g^*$ is said to be \emph{$\eta$-polarized} 
if its $\eta$th component, $\Psi^\eta \colon M \to \R$, is proper 
and bounded from below.
\end{Definition}

Note that an $\eta$-polarized function is necessarily proper,
because its $\eta$-component is proper.
However, not every proper map is $\eta$-polarized for some $\eta$.
If $M$ is compact, $\Psi$ is automatically $\eta$-polarized for all $\eta$.

Polarized abstract moment maps posses the following two properties,
which may in generally fail for proper abstract moment maps:
\begin{enumerate}
\item
A linear combination of $\eta$-polarized abstract moment maps on the same 
manifold is again an $\eta$-polarized (hence proper) abstract moment map.
\item
Let $\Psi_j\colon M_j \to \g^*$, $j=1,2$, be $\eta$-polarized abstract moment maps.
Consider the product $M_1 \times M_2$ with the diagonal $G$-action;
let $\pi_1, \pi_2$ be the projection maps to $M_1$ and $M_2$. Then
$\Psi_1 \circ \pi_1 + \Psi_2 \circ \pi_2$ is an $\eta$-polarized (hence proper)
abstract moment map.
\end{enumerate}

Fix a $G$-manifold $M$ and a vector $\eta \in \g$. 
Denote the set of zeros of $\eta_M$ by $M^\eta$.
This is exactly the set of points whose infinitesimal stabilizer contains 
$\eta$.  Therefore, the $\eta$-coordinate of any assignment is well defined 
on this set.

\begin{Definition}
An assignment $A$ is \emph{$\eta$-polarized} if its $\eta$-component
on $M^\eta$,
$$ A^\eta\colon M^\eta \to \R,$$
is proper and bounded from below.
\end{Definition}

\comment{\ Omit the following remark. I'm keeping it in $\backslash$comment
in case we want to use it for the book.\\ 
\begin{Remark}
Let $H$ be the closure in $G$ of the one-parameter subgroup generated by
$\eta$. Then the zero set $M^\eta$ of $\eta_M$ is equal to the fixed point 
set $M^H$ of $H$ in $M$. On this set, the $H$-component of an assignment 
is a well--defined function $A^H \colon  M^H \to \h^*$. 
An assignment $A$ is $\eta$-polarized if and only if this function
is $\eta$-polarized in the sense of Definition \ref{polarized}. 
For a generic $\eta$, we have $H=G$. For such an $\eta$, an assignment is 
$\eta$-polarized if its restriction to the fixed point set, which is
an ordinary function $A\colon M^G \to \g^*$, is $\eta$-polarized
in the sense of Definition \ref{polarized}.
\end{Remark} }

\begin{Theorem} \labell{thm:assignment proper}
Let $M$ be a manifold with a $G$ action.
For every $\eta$-polarized assignment $X \mapsto A(X)$ on $M$ 
there exists an $\eta$-polarized abstract 
moment map $\Psi \colon M \to \g^*$ whose assignment is $A$.
\end{Theorem}

\begin{Corollary} \labell{cor:torus-exist}
Assume that $M^G$ is compact. Then every assignment extends to a proper
abstract moment map.
\end{Corollary}

As a consequence, if $M^G$ is compact, there always exists a proper
abstract moment map (e.g., one which extends the zero assignment).

\begin{proof}[Proof of Theorem \ref{thm:assignment proper}]
The function $A^\eta\colon M^\eta \to \R$ is well defined, proper, and
bounded from below. Since $M^\eta$ is closed, $A^\eta$ extends to
a function $\varphi \colon M \to \R$ that is proper and bounded
from below. (See item 1 in the proof of Theorem \ref{thm:circle-exist}.)

For each $m \in M$, let $\Psi_m \in \g^*$ be an element whose projection 
to $\g_m^*$ is $A(m)$.
We choose 
$\Psi_m \in \g^*$ to meet the following additional requirement: 
$\Psi_m^\eta=\left< \Psi_m , \eta \right> = \varphi(m)$. 
If $m \in M^\eta$, this condition is automatically satisfied,
and if $m \not \in M^\eta$, 
this choice is possible because $\eta\not\in\g_m$.

Let $U_m$ be a tubular neighborhood of
the orbit through $m$ which equivariantly retracts to the orbit and 
on which the function $\varphi$ differs from the value $\varphi(m)$ 
by less than $1$.
Then the constant function $\Psi_m$ is an abstract moment map on $U_m$ with 
assignment $A$ and whose $\eta$-component is bounded from below
by $\varphi-1$. 

Choose an invariant partition of unity $\{\rho_j\}$ subordinate to the
covering of $M$ by the open subsets $U_m$, with the support of $\rho_j$
contained in the open set $U_{m_j}$.
Then the convex combination
$\Psi = \sum \rho_j \Psi_{m_j}$ is an abstract moment map;
this follows from Lemma \ref{fix A}, applied to open subsets of the
manifold. Moreover, since the $\eta$-component of each $\Psi_m$
is bounded from below by $\varphi-1$, the same holds for $\Psi$.
Since $\Psi^\eta \geq \varphi-1$, and 
$\varphi-1$ is proper and bounded from below,
$\Psi^\eta$ is proper and bounded from below. 
\end{proof}

In general, a $G$-manifold $M$ may admit no proper abstract moment maps, 
even when every connected component of the fixed point set $M^G$ is compact 
(and so a proper locally constant map $\psi\colon M^G\to\g^*$ does exist).
The obstruction lies in the compatibility condition; the manifold $M$ 
might not admit a proper assignment. 
We will now construct an example of such a $G$-manifold.

\begin{Example}
\labell{S4-chain}
Let $G = S^1 \times S^1$ act on the four--dimensional sphere $S^4$
as in Example \ref{S4}. Recall that the fixed points are the North 
and South Poles and that any abstract moment map on $S^4$
must take the same value at these points.
Fix some small  $\epsilon>0$, and let  $D^4$  be the $\epsilon$-ball
in  $\C\times\C$ with the $G$-action that rotates each of the two factors.
Take the trivial disk bundle over $S^4$,
\begin{eqnarray*}
 N & = & S^4\times D^4  \\
   & = & \{ (z,z',x,w,w') \ | \ |z|^2 + |z'|^2 + x^2 =1 \text{ and }
                                |w|^2 + |w'|^2 < \epsilon^2 \} \\
   & \subset & \C^2 \times \R \times \C^2
\end{eqnarray*}
with the diagonal action of $G$.
Since the neighborhood of each of the two fixed points in $N$ is
equivariantly diffeomorphic to  $D^4\times D^4$,  we can plumb an 
infinite sequence of such  $N$'s.  More explicitly, take  
$M= N\times\Z/\!\sim$  where the equivalence  relation $\sim$ is
\begin{equation} \labell{gluing-map}
 (z,z',x,w,w',n) \sim (w,w',-x,z,z',n+1) 
\end{equation}
for all $x>0$ and $n \in \N$, whenever both $|z|^2+|z'|^2$ and 
$|w|^2+|w'|^2$ are less than $\epsilon$.  Then $M$ is a $G$-manifold.
The gluing map \eqref{gluing-map} reverses the orientation;
however, we can get an orientation on $M$ by flipping the
orientation of every other copy of $N$.
An abstract moment map on $M$ must take a constant value on the infinite
sequence of fixed points; such a map cannot be proper.
\end{Example}

\subsection{Minimal stratum assignments.}
\labell{sec:minimal} 
Theorem \ref{thm:assignment} can be understood as that assignments
are combinatorial counterparts of abstract moment maps. The amount
of information needed to determine an assignment can be further
reduced by taking a full advantage of the compatibility condition,
as follows. Recall that the (infinitesimal) orbit type strata in $M$ 
are partially ordered;
$X \preceq Y$ if and only if $X$ is contained in the closure of $Y$.
The strata that are minimal under this ordering
are exactly those that are closed subsets of $M$.
The closure of any orbit type stratum in $M$ is a smooth sub-manifold which
contains a minimal stratum.

Every component of the fixed point set, $M^G$, is a minimal stratum.
However, there can exist minimal strata outside the fixed point set $M^G$.  
Whether or not such strata exist
is related to an algebraic property called \emph{formality}.
Recall that a compact manifold $M$ is formal if one of the following 
equivalent conditions (see, e.g., \cite{Borel}, \cite{hsiang}, or 
\cite{Ki:book}) is satisfied:
\begin{enumerate}
\item $H^*_G(M)=H^*(M)\otimes H^*(BG)$ as an $H^*(BG)$-module;
\item $H^*_G(M)$ has no $H^*(BG)$-torsion; 
\item the restriction 
$j^*\colon H^*_G(M)\to H^*_G(M^G)=H^*(M^G)\otimes H^*(BG)$
is a mono-morphism.
\end{enumerate}

Here is an interesting geometric consequence of formality:

\begin{Proposition}
\labell{exam:min-strata-formal}
On a compact formal manifold $M$, every minimal stratum is a connected 
component of $M^G$. 
\end{Proposition}

\begin{proof}
Let $X$ be a minimal stratum and let $H$ be the connected component
of identity of $G_x$ for $x\in X$.  Assume $H \neq G$.
Then the equivariant Thom class $\tau$ of the normal bundle to $X$ is a 
non-zero torsion element in $H^*(BG)$. In fact, $\tau$ is
annihilated by the image of $H^*(B (G/H))\to H^*(BG)$. 
Alternatively, $j^*\tau=0$, because $X\cap M^G=\emptyset$.
\end{proof}

For example, when $M$ is compact symplectic with $G$ acting
Hamiltonianly, $M$ is equivariantly perfect, and hence formal,
(see \cite{Ki:book}), and 
the above analysis applies. In this case, however, to show that 
a minimal stratum $X$ consists of fixed points,
it suffices to observe that $X$ is a compact 
symplectic manifold and $H$ acts Hamiltonianly on $X$, so $H$ must have
fixed points on $X$.

\comment{\ Remove the following Question. (Perhaps keep it as
$\backslash$comment for our own use.)\\ 
\begin{Question}
It would be interesting to find out to what extend the equivariant cohomology
$H_G^*(M)$ can be recovered from the set of orbit type strata with its
partial ordering and with their
stabilizers or, alternatively, from the
\emph{isotropy assignment} of Section \ref{subsec:other coefficients}.\\ 
The next question is a bit problematic: spectral sequences, etc.
We need to think more of it before we ask it.\\ 
Moreover, can we reconstruct $H_G^*(M)$ if we know $H_G^*(X)$ 
for all the strata $X$, as well as the partial ordering between strata?
\end{Question} 
Compare with what Kirwan does in her book. Also, compare with
Bredon's theorem. He says something like that the equivariant homotopy
type of a manifold is determined by the ordinary homotopy types
of the strata and by their partial ordering.}

\begin{Definition} 
A \emph{minimal stratum assignment} 
is an assignment of an element $A(X)\in\g_X^*$ to each minimal 
stratum $X$, where $\g_X$ is the infinitesimal
stabilizer of $x\in X$, such that the following compatibility condition
is satisfied:
\emph{if two minimal strata $X_1$ and $X_2$ are such that
$X_1\preceq Y$ and $X_2\preceq Y$ for some stratum $Y$, 
then the restrictions to $\g_Y$
of $A(X_1)$ and of $A(X_2)$ are the same: 
$A(X_1)^{\g_Y}=A(X_2)^{\g_Y}$.}
\end{Definition}

Notice that this condition holds automatically for the zero assignment.

The following theorem follows immediately from the definitions.

\begin{Theorem} \labell{min-assign}
The restriction of any assignment to the minimal strata
is a minimal stratum assignment. Conversely, any minimal stratum 
assignment extends to a unique assignment. Hence,
every minimal stratum assignment is associated with an abstract moment map.
\end{Theorem}

\begin{Remark}
It appears that in Theorem \ref{min-assign} the minimal stratum assignment 
cannot be replaced by a function defined only on the fixed point set.
Namely, we expect there to exist a $G$-manifold $M$ with isolated 
fixed points and a function $\psi \colon M^G \to \g^*$ which does not
extend to an assignment (hence does not extend to an abstract moment map),
but which satisfies the following compatibility condition:
if $x,y \in M^G$ belong to the same connected component 
of $M^H$, then $\psi^H(x) = \psi^H(y)$.
\end{Remark}

\begin{Question}
In Remark \ref{rmk:non-abelian} we proposed a definition of
assignments for an action of not necessarily abelian Lie group.
It appears to be an interesting and feasible problem to check whether
or not the results of this section generalize to such actions.
\end{Question}

\section{Exact moment maps}
\labell{sec:exact} 

We have already shown that the natural forgetful
homomorphism from the space of abstract moment maps on a 
$G$-manifold $M$ to the space $\AA(M)$ of assignments on $M$
is onto (Theorem \ref{thm:assignment}). In this section
we study the kernel of this epimorphism.

Recall that an abstract moment map $\Psi$ is said to be exact
if there exists a $G$-invariant one-form $\mu$ with 
$\Psi^\xi=\mu(\xi_M)$ for all $\xi\in\g$ (see Example \ref{ex:one-form}).
The assignment associated with such a map is zero.
The following result, which is proved later in this section,
shows that the converse is also true.

\begin{Theorem} \labell{cor1}
An abstract moment map whose assignment is identically zero is exact.
More explicitly, suppose that $\Psi\colon M \to \g^*$ 
is an abstract moment map such that for each subgroup 
$H \subset G$, the function $\Psi^H \colon M \to \h^*$ vanishes on the 
$H$-fixed point set $M^H$. Then there exists an invariant one-form $\mu$
such that $\Psi^\xi = \mu(\xi_M)$ for all $\xi \in \g$.
\end{Theorem}

Combining Theorems \ref{thm:assignment} and \ref{cor1} we obtain

\begin{Corollary}
\labell{cor12}
The sequence
$$
0\to
\left\{ \begin{array}{c}
\text{exact}\\ \text{moment} \\ \text{maps} 
\end{array} \right\}
\to
\left\{ \begin{array}{c}
\text{abstract} \\ \text{moment} \\ \text{maps}
\end{array} \right\}
\to \AA(M)\to 0
$$
is exact.
\end{Corollary}

The proof of Theorem \ref{cor1} relies on the following key result,
which we will prove in Section \ref{sec:proof}:

\begin{Theorem}
\labell{twoform-linear}
Let $G$ be a torus acting linearly on $\reals^m$,
and let $\Psi$ be an abstract moment map on a neighborhood of the origin,
vanishing at the origin.
Then there exists a $G$-invariant one-form 
$\mu$ on a neighborhood of the origin such that $\mu(\xi_M)=\Psi^\xi$
for all $\xi \in \g$.
\end{Theorem}

We will also need a parametric version of this theorem:
\begin{Corollary}
\labell{twoform-linear-par}
Let $G$ be a torus acting linearly on the fibers of a vector bundle 
$\Nu \to Y$, and let $\Psi$ be an abstract moment map on a neighborhood 
of the zero section, vanishing on the zero section. Then there exists 
a smooth family $\mu$ of $G$-invariant one-forms on the fibers of $\Nu$,
such that $\mu(\xi_M) = \Psi^\xi$ near the zero section.
\end{Corollary}

\begin{proof}[Proof of Corollary \ref{twoform-linear-par}]
By using a partition of unity on $Y$, the corollary can be reduced 
to the case where $Y$ is a linear space and $\Nu = \R^m \times Y$.
This case follows immediately from Theorem \ref{twoform-linear} 
when $\R^m$ is replaced by $\R^m\times Y$ 
with the trivial $G$-action on the second factor.
\end{proof}

Assuming Theorem \ref{twoform-linear}, let us prove a preliminary result,
which is a local version of Theorem \ref{cor1} that will be used in the
next section, and deduce Theorem \ref{cor1} from it.

\begin{Proposition} 
\labell{prop:one-form}
Let $G$ be a torus acting on a manifold $M$ and let $\Psi \colon M \to \g^*$
be an abstract moment map. Let $p$ be a point in $M$ and $H = G_p$ 
its stabilizer.  Suppose that $\Psi^H(p)=0$. 
Then there exists an open $G$-invariant neighborhood $V$ of $p$ in $M$ 
and a $G$-invariant one-form $\mu$ on $V$ such that 
\begin{equation}
\labell{eq:one-form}
\mu(\xi_M) = \Psi^\xi
\end{equation}
on $V$ for all $\xi \in \g$.
\end{Proposition}

\begin{proof}
Let us first examine the case where the action is locally free near $p$. 
Fix a basis $\xi_1,\ldots,\xi_n$ in $\g$. Then the vector fields
$(\xi_i)_M$ form a basis in the tangent space to the orbit at every point
of a $G$-invariant neighborhood $V$ of the orbit through $p$. By setting 
\begin{equation}
\labell{eq:loc-free}
\mu((\xi_i)_M)=\Psi^{\xi_i} , 
\end{equation}
we thus obtain a form
defined along the orbits in $V$. We extend it to a differential form
$\mu$ on $V$ by taking its composition with an orthogonal projection to
the orbit with respect to a $G$-invariant metric. It is easy to see
that $\mu$ satisfies the condition $\mu(\xi_M) =\Psi^\xi$.

Let us now prove the proposition in the general case.
Pick a closed subgroup $K \subset G$ whose Lie algebra $\k$ is 
complementary to $\h$ in $\g$.  
A small $G$-invariant neighborhood $V$ of the orbit $Y$ through $p$
can be identified, by the slice theorem, with a neighborhood of the
zero section in the normal bundle $ \pi \colon \Nu \to Y$ to $Y$ in $M$,
with the action induced by that on $M$. We can apply Corollary 
\ref{twoform-linear-par} to the linear $H$-action on the fibers of $\Nu$,
equipped with the abstract moment map $\Psi^H$ induced from $M$.
As a result, we get a smooth family $\mu$ of one-forms 
on the fibers of $\Nu$, such that $\mu(\xi_M) = \Psi^\xi$ for all $\xi\in\h$.
The $K$-orbits form a foliation which is transverse to the fibration $\pi$.
We extend $\mu$ to a one-form on a whole neighborhood of $Y$
by making $\mu$ vanish on the vectors tangent to the $K$ orbits.
The resulting form is a $G$-invariant form $\mu_H$ on $V$ so that 
$\mu_H(\xi_M) = \Psi^\xi$ for all $\xi \in \h$, and
$\mu_H(\xi_M) = 0$ for all $\xi \in \k$.

The $K$-action on $V$ is locally free. 
Let $\mu_K$ be the form defined as above by \eqref{eq:loc-free} 
and extended to $V$ so that it vanishes on the vectors tangent
to the fibers of $\pi$. Then 
$\mu_K(\xi_M) = \Psi^\xi$ for all $\xi \in \k$.
Since the vector fields $\xi_M$ for $\xi\in\h$ are tangent
to the fibers of $\pi$, we also have 
$\mu_K(\xi_M) = 0$ for all $\xi \in \h$.
The form $\mu=\mu_H+\mu_K$ has the desired property, 
that $\mu(\xi_M) = \Psi^\xi$ for all $\xi \in \g$.
\end{proof}

\begin{proof}[Proof of Theorem \ref{cor1}]
By Proposition \ref{prop:one-form}, there exists an open covering
of $M$ by invariant sets $U_\alpha$, and on each $U_\alpha$ there
exists an invariant one-form $\mu_\alpha$ such that 
$\Psi^\xi = \mu_\alpha(\xi_M)$ for all $\xi \in \g$. Let $\rho_j$
be a partition of unity subordinate to this covering, with $\rho_j$
supported in $U_{\alpha_j}$ for each $j$. Define 
$\mu = \sum \rho_j \mu_{\alpha_j}$. Then $\Psi^\xi = \mu(\xi_M)$
on $M$ for all $\xi \in \g$.
\end{proof}

\section{Hamiltonian moment maps}
\labell{sec:forms} 

We have already seen (Example \ref{exam:presymplectic}) that 
every moment map which is associated with a closed invariant two-form
is an abstract moment map.
We will examine now the question of which abstract
moment maps arise in this way. Recall that such abstract moment 
maps are called Hamiltonian. Thus we fix an abstract moment map $\Psi$,
and we look for a closed two-form $\omega$ with
$\iota(\xi_M) \omega = - d \Psi^\xi $. 
Note that such an $\omega$ would necessarily be $G$-invariant.

Our first observation is an immediate consequence of the fact
that every exact moment map, $\Psi^\xi=\mu(\xi_M)$, is automatically 
Hamiltonian with $\omega=d\mu$.  Thus Theorem \ref{cor1} implies

\begin{Corollary} \labell{cor11}
Let $\Psi \colon M \to \g^*$ be an abstract moment map with zero assignment.
Then $\Psi$ is associated with an exact two-form. 
In particular, $\Psi$ is Hamiltonian.
\end{Corollary}

\subsection{Local existence of two-forms.}
Our next result shows that there are no local obstructions to
the existence of $\omega$, if $G$ is abelian. 
Quite surprisingly, a similar local existence result fails to hold 
for non-abelian compact groups, \cite{Br}.

\begin{Corollary}[Local existence of two-forms]
\labell{thm:presymplectic}
Let $G$ be a torus acting on a manifold $M$,
and let $\Psi \colon M \to \g^*$ be an abstract moment map.
For every $p\in M$, $\Psi$ is associated with an exact two-form $\omega$ 
on some open $G$-invariant neighborhood $V$ of $p$.
In particular, $\Psi$ is Hamiltonian on a neighborhood of $p$.
\end{Corollary}

\begin{proof}
Consider the new abstract moment map $\Psi - \Psi(p)$.
By Proposition \ref{prop:one-form}, there exists an invariant
neighborhood $V$ of $p$ in $M$ and a $G$-invariant one-form $\mu$
on $V$ such that $\mu(\xi_M) = \Psi^\xi - \Psi^\xi(p)$ on $V$.
Let $\omega = d\mu$; then $\iota(\xi_M)\omega = -d\Psi^\xi$ on $V$.
\end{proof}

The following semi-local result is also of interest.

\begin{Corollary} \labell{cor2}
On a manifold with a unique minimal stratum, $X$,
every abstract moment map is Hamiltonian.
\end{Corollary}

\begin{proof}
Pick an element $\gamma \in \g^*$  whose restriction to $\g_X$
is equal to $\Psi^{\g_X}(X)$, and apply Theorem \ref{cor1}
to the abstract moment map $\Psi - \gamma$.
\end{proof}

\comment{\ Ax in the paper and keep for the book:\\ 
\begin{Corollary}
In every $G$-manifold, there exists an invariant neighborhood 
of the fixed point set, $M^G$, on which every abstract moment map
is Hamiltonian.
\end{Corollary} }

\begin{proof}
Take a tubular neighborhood of $M^G$ which retracts to $M^G$, and apply 
Corollary \ref{cor2} to each of its connected components.
\end{proof}

\begin{Remark}
It is well known that moment maps $\Psi$ associated with symplectic forms
satisfy a certain non-degeneracy condition. For example, for
circle actions the Hessian $d^2\Psi$ must be non-degenerate on
the normal bundle to the fixed point set. In \cite{GGK:book}, we 
state explicitly a necessary and sufficient condition for $\Psi$
to be locally, near $M^G$, associated with a 
symplectic form. Furthermore, we will prove that abstract moment maps
satisfying this non-degeneracy condition globally have many properties
of moment maps on symplectic manifolds. These include
the convexity theorem (\cite{At:convexity} and \cite{GS:convexity})
and formality (\cite{Ki:book}, see also Section \ref{sec:minimal} 
above). 
\end{Remark}

\subsection{Global existence of two-forms}
\labell{sec:forms-global}
Let us now turn to the problem of global existence for $\omega$. 
The following example shows that not every abstract moment map
is Hamiltonian.

\begin{Example} \labell{CP2}
Let $S^1$ act on $\CP^2$ by 
$$
\lambda \cdot [z_0:z_1:z_2] = [z_0: \lambda z_1: \lambda^2 z_2].
$$
There are three fixed points: $[1:0:0]$, $[0:1:0]$, and $[0:0:1]$.
Denote by $a,b,c$ their respective images by an abstract moment map.
If the abstract moment map is associated with a closed two form, 
$\omega$, then it is an easy consequence of Stokes's theorem 
that the differences,
$b-a$ and $c-b$ are, respectively, equal (up to a factor) 
to the integrals of $\omega$
on the 2-spheres $[*:*:0]$ and $[0:*:*]$ in $\CP^2$. Since these
lie in the same cohomology class, the values $a,b,c$ must then be
equidistant: $a-b=b-c$. However, an abstract moment map can take 
arbitrary values $a,b,c$ at the three fixed points, by 
Theorem \ref{thm:circle-exist}.
\end{Example}

Recall that the equivariant cohomology classes in $H^2_G(M)$
are represented by the sums $\omega + \Psi$
where $\Psi$ is a Hamiltonian moment map and $\omega$ is a 
compatible two-form; see Example \ref{exam:presymplectic}. 
The forgetful mapping which sends $\omega+\Psi$ to the assignment 
$A$ corresponding to $\Psi$ gives rise to a homomorphism 
$$\rho\colon H^2_G(M)\to  \AA(M).$$

\begin{Theorem}
\labell{thm:exist-2form}
An abstract moment map $\Psi$ is Hamiltonian if and only if 
$A\in \im \rho$, where $A$ is the assignment of $\Psi$. 
\end{Theorem}

\begin{proof}
It is clear by definition that $A\in\im\rho$ if $\Psi$ is Hamiltonian.

Conversely, assume that $A\in \im\rho$. Then there exists
a $G$-equivariant equivariantly closed two-form $\omega+\Phi$ such
that the assignment of $\Phi$ is also $A$.
The difference $F=\Psi-\Phi$ is an abstract moment map with the zero
assignment. By Theorem \ref{cor1}, $F$ is exact and therefore
Hamiltonian (Corollary \ref{cor11}). Thus $\Psi$ is Hamiltonian as the
sum of two Hamiltonian abstract moment maps, $F$ and $\Phi$.
\end{proof}

The space of Hamiltonian assignments, i.e., 
assignments associated with Hamiltonian abstract moment maps, 
is the quotient of the space of all Hamiltonian abstract moment
maps by the space of exact abstract moment maps.
This follows from Corollary \ref{cor12}.
These three spaces fit together to form a part of a commutative square of
exact sequences which summarizes some of our results.

\begin{Proposition}
\labell{prop:diagram}
The following diagram is commutative and all of its rows and columns
are exact:
$$\begin{array}{ccccccccc}
&  & 0 &  & 0 &  & 0 &  & \\
&  & \downarrow &  & \downarrow &  & \downarrow &  & \\ 
0&\to &
\left\{ \begin{array}{c}
\text{basic}\\ \text{exact}\\ \text{2-forms}
\end{array}\right\}&
\to &
\left\{\begin{array}{c}
\mbox{equivariantly}\\ \mbox{exact}\\ \mbox{2-forms}
\end{array}\right\}&
\to &
\left\{ \begin{array}{c}
\mbox{exact} \\ \mbox{moment} \\ \mbox{maps}
\end{array} \right\}&
\to & 0\\
  &  & \downarrow &  & \downarrow &  & \downarrow &  &  \\ 
0&\to &
\left\{ \begin{array}{c}
\text{basic}\\ \text{closed}\\ \text{2-forms}\end{array}\right\}&\to &
\left\{\begin{array}{c}
\mbox{equivariantly}\\ \mbox{closed}\\ \mbox{2-forms}
\end{array}\right\}&
\to &
\left\{ \begin{array}{c}
\mbox{Hamiltonian} \\ \mbox{moment} \\ \mbox{maps}
\end{array} \right\}&
\to & 0\\
  &  & \downarrow &  & \downarrow &  & \downarrow &  &  \\ 
0 & \to & H^2(M/G) & \to & H^2_G(M) & \to &
\left\{ \begin{array}{c}
\mbox{Hamiltonian} \\ \mbox{assignments}
\end{array}\right\}& \to & 0 \\
  &  & \downarrow &  & \downarrow &  & \downarrow &  &  \\ 
&  & 0 &  & 0 &  & 0 &  & 
\end{array}
$$
\end{Proposition}

\begin{proof}
The exactness of the left column is a particular case of a more general 
fact, that the cohomology of the basic De Rham complex of $M$
is equal to $H^*(M/G)$, if $G$ is compact or, more generally, 
if the $G$-action is proper, even when the action is not free. 
This result, due to Koszul \cite{Koszul}, is similar to the De Rham theorem 
and can be proved in the same way. An easy proof is as follows. 
Recall that the sequence of
sheaves of singular cochains on $M/G$ (with real coefficients) 
is a fine resolution of the constant sheaf $\R$ on $M/G$.
Furthermore, basic forms on $G$-invariant open subsets of $M$ form a sheaf
on $M/G$. This sheaf is a resolution of the locally constant sheaf because
it is locally acyclic. Indeed, by using the fact that $G$ is compact 
(or that the action is proper) and adapting the proof of the Poincar\'{e}
lemma, one can show that the basic cohomology of a neighborhood of an
orbit is the same as of the orbit itself, i.e., zero in positive degrees.
It is easy to see that this sheaf is also fine because it admits partitions
of unity. Thus basic forms on $M$ provide another fine resolution of 
the constant sheaf on $M/G$.
Since the cohomology of both resolutions
are equal to the \v{C}ech cohomology of the constant sheaf on $M/G$,
they are equal to each other.

The middle column is exact by the definition of equivariant cohomology
via the equivariant De Rham complex. 
(See, e.g., \cite{AB} and \cite{DKV}.)

Exactness of the right column follows from Corollary \ref{cor12}.

The fact that the top two rows are exact follows directly from the 
definitions of the spaces involved.

The commutativity of the diagram is clear. Finally, commutativity
with the exactness of the columns and the top two rows 
implies that the bottom row is exact by simple diagram chasing.
\end{proof}

\begin{Remark}
Our notion of assignments has an interesting connection with 
a recent theorem of Goretsky-Kottwitz-MacPherson, \cite{GKM}.
Assume that a compact oriented $G$-manifold $M$ is formal and satisfies
in addition the so-called GKM condition, which we will recall below.
Then
the Goretsky-Kottwitz-MacPherson theorem implies that every assignment 
is Hamiltonian,
and hence every abstract moment map is associated with a two-form.

The ``GKM condition" is that the fixed points are isolated, 
and, additionally, every orbit type stratum with stabilizer 
of codimension one is two-dimensional.

Let us now show how the above assertion follows from the theorem.
Let $A$ be an assignment.  Its restriction to the fixed point set
is a locally constant function, $A^G\colon M^G \to \g^*$. Such a function
can be identified with an element of $H^2_G(M^G)$. For each subgroup 
$H \subset G$ of codimension one, the set of $H$-fixed points is a disjoint
union of two-spheres in $M$, on each of which $G/H$ acts with exactly
two fixed points; this is a consequence of the ``GKM condition"
and formality. The assignment compatibility condition implies that
in each such a two-sphere, the images in $\h^*$ of $A^G$ are the same
at the two fixed points. The Goretsky-Kottwitz-MacPherson theorem
asserts that this condition on $A^G$ implies that there exists
an equivariantly closed equivariant two-form, $\omega+\Psi$, on $M$,
whose restriction to $M^G$ is $A^G$. By formality, an assignment
on $M$ is uniquely determined by its restriction to $M^G$
(see Proposition \ref{exam:min-strata-formal}).
Hence, $A$ is the assignment associated with $\Psi$; 
hence, it is Hamiltonian.
\end{Remark}

\section{Abstract moment maps on linear spaces: \\ the proof of Theorem
\ref{twoform-linear}.}
\labell{sec:proof}

A crucial step in the proof of Theorem \ref{cor1}
and Proposition \ref{prop:one-form}, on  which many of our
subsequent results rely, is Theorem \ref{twoform-linear}. In this
section we recall this theorem and prove it. A different proof
can be found in \cite{GGK:book}. 

\smallskip \noindent \textbf{Theorem \ref{twoform-linear}.}
\emph{
Let $G$ be a torus acting linearly on $\reals^m$,
and let $\Psi$ be an abstract moment map on a neighborhood of the origin,
vanishing at the origin.
Then there exists a $G$-invariant one-form 
$\mu$ on a neighborhood of the origin such that $\mu(\xi_M)=\Psi^\xi$
for all $\xi \in \g$.} 

\begin{proof}[Proof of Theorem \ref{twoform-linear}]
First note that by adding, if necessary, an additional copy of $\R$
(with the trivial $G$-action) to $\R^m$, we can always make
$m$ even. 

There exists a $G$-invariant complex structure on the vector space $\R^m$;
fix one.
We obtain a representation $\rho$ of $G$ on $\C^d$
with weights $\alpha_1,\ldots,\alpha_d$.
The infinitesimal action of the Lie algebra $\g$ of $G$ on $\C^d$ is
then given by the vector fields
\begin{equation}
\labell{eq:1}
\xi_M = \sqrt{-1} \sum \limits_{i=1}^d \alpha_i(\xi) 
 \left(z_i\frac{\p}{\p z_i}-\bar z_i\frac{\p}{\p\bar z_i}\right) .
\end{equation}

 From now on we will forget about the reality conditions and assume 
that $\mu$ is a complex-valued one-form. When such a form is found, it
will suffice to replace it by the real form $(\mu+\bar{\mu})/2$, 
which still has the desired properties 
because $\xi_M$ are real vector fields.

We will also forget about the equivariance conditions. 
If a form $\mu$ such that 
$\mu(\xi_M) = \Psi^\xi$ for all $\xi \in \g$ is constructed,
and if $\Psi$ is equivariant, we can replace $\mu$ by its average.

Any one-form on $\C^d$ can be written in the form
\begin{equation} \labell{mu on Cd}
\mu=-\sqrt{-1} \sum\limits_{i=1}^d f_i dz_i-g_i d\bar z_i
\end{equation}
for some smooth functions $f_j$ and $g_j$, $j=1,\ldots,d$.
For such a one-form, the function $\Psi \colon \C^d \to \g^*$
defined by $\Psi^\xi = \mu(\xi_M)$,
with the vector fields $\xi_M$ given by \eqref{eq:1}, is
\begin{equation} \labell{Psi on Cd}
\Psi = \sum_{j=1}^d (z_j f_j + \zbar_j g_j) \alpha_j.
\end{equation}
Conversely, for any function $\Psi$ which has the form 
\eqref{Psi on Cd} for some smooth functions $f_j$, $g_j$, there exists
a one-form $\mu$ such that $\Psi^\xi = \mu(\xi_M)$ for all $\xi \in \g$.
Namely, just take $\mu$ to be given by \eqref{mu on Cd}.

To prove Theorem \ref{twoform-linear}, it is thus enough to prove
the following
\begin{Proposition} \labell{smooth}
Let $\Psi$ be a $\g^*$-valued function on a neighborhood of the origin
in $\C^d$, vanishing at the origin and satisfying the second condition of 
an abstract moment map: 
\begin{quote}
for any subgroup $H \subset G$, 
the function $\Psi^H \colon \C^d \to \h^*$ is locally constant on the set
of $H$-fixed points.
\end{quote}
Then there exist smooth functions $f_j$ and $g_j$
such that $\Psi$ is given by \eqref{Psi on Cd} on a neighborhood of the
origin.
\end{Proposition}

\begin{Remark}
The converse is easy: any $\Psi$ of the form \eqref{Psi on Cd}
satisfies the second condition of an abstract moment map.
\end{Remark}

Let us start with polynomial functions and one-forms:

\begin{Proposition} \labell{polynomial}
Let $\Psi \colon \C^d \to \g^*$ be a polynomial function which vanishes
at the origin and which satisfies
the second condition of an abstract moment map. Then there exist
polynomials $f_j$ and $g_j$ on $\C^d$ such that $\Psi$ is given by
\eqref{Psi on Cd}.
\end{Proposition}

\begin{proof}[Proof of Proposition \ref{polynomial}.]
Since $\Psi$ is polynomial, we can write it uniquely as a sum of monomials,
$$ \Psi = \sum\limits_{k,l} \beta_{k,l} z^k \zbar^l,$$
summing over $k = (k_1 , \ldots, k_d)$ and $l = (l_1,\ldots,l_d)$
in $\N^d$, where the coefficients $\beta_{k,l}$ are in $\g^*$.

For every subset $I \subset \{ 1, \ldots, d\}$, denote by $(\Cx)^I$
the subset of $\C^d$ consisting of all vectors $(z_1,\ldots,z_d)$
for which $z_i \neq 0$ if and only if $i \in I$.
All the points $z$ in $(\Cx)^I$ have the same stabilizer, $G_I$,
whose Lie algebra is 
\begin{equation} \labell{gI}
\g_I = \bigcap\limits_{i \in I} \ker \alpha_i.
\end{equation}
Since $\Psi$ satisfies the second condition of an abstract moment map,
$\Psi^{\g_I}$ is constant on $(\Cx)^I$. Since, additionally, $\Psi$
is continuous and vanishes at the origin,
$\Psi^{\g_I}$ vanishes on $(\Cx)^I$.

Let us analyze what this condition tells us about the coefficients
$\beta_{k,l}$. The polynomial
$\Psi^\xi = \sum_{k,l} \beta_{k,l}(\xi) z^k \zbar^l$ 
vanishes on $(\Cx)^I$ for all $\xi \in \g_I$
if and only if for each $k,l$, the summand
$\beta_{k,l}(\xi) z^k \zbar^l$ vanishes on $(\Cx)^I$
for all $\xi \in \g_I$.  Fix $k$ and $l$, and restrict attention to
$$I = I_{k,l} = \{ i \ | \ k_i \neq 0 \text{ or } l_i \neq 0 \}.$$
Since the monomial $z^k \zbar^l$ does not vanish on $(\Cx)^I$,
its coefficient, $\beta_{k,l} (\xi)$, must vanish for all $\xi \in \g_I$.
By \eqref{gI}, a linear functional that vanishes on $\g_I$ 
is a linear combination of $\alpha_i$, $i \in I$. Therefore,
$\beta_{k,l} = \sum\limits_{i \in I_{k,l}} \lambda_{i,k,l} \alpha_i$,
and
$$ \Psi = \sum_i \alpha_i
   \sum\limits_{k,l \text{ such that } i \in I_{k,l} }
   \lambda_{i,k,l} z^k \zbar^l.$$
Since for each $i \in I_{k,l}$, either $z_i$ or $\zbar_i$ factors out of
the monomial $z^k \zbar^l$, $\Psi$ is of the form \eqref{Psi on Cd}.
\end{proof}

We will now show that the theorem we need to prove in the smooth category
follows from its polynomial version, which has already been proved.
In other words, we will deduce Proposition \ref{smooth} from Proposition 
\ref{polynomial}. To this end, let us reformulate these propositions 
as assertions that certain sequences of homomorphisms are exact.

Denote by $\PP$ the ring of complex-valued polynomials in $z_j$ 
and $\bar z_j$, $j=1,\ldots,d$.
Define the modules $\MM_i$, $i=1,\ 2, \ 3$, over $\PP$ as follows:
\begin{itemize}
\item $\MM_1$ is the space of one-forms 
$\sum f_i\,dz_i+g_i\,d{\bar z}_i$ with $f_i$ and $g_i$ in $\PP$.

\item $\MM_2$ is the tensor product $\PP \otimes \g^*$ over $\C$.

\item For each subset $I \subseteq \{ 1, \ldots, n \}$,
denote by $\PP_I$ the ring of polynomial functions on $(\Cx)^I$.
The restriction homomorphism $\PP \to \PP_I$ makes $\PP_I$ into a 
$\PP$-module. Set
$$
\MM_3 = \bigoplus_I \PP_I \otimes \g_I^* ,
$$
where $\g_I$ is the Lie algebra of the stabilizer of $(\Cx)^I$,
given by \eqref{gI}.
\end{itemize}
Define the sequence of homomorphisms
\begin{equation}
\labell{eq:(2)}
\MM_1\overset{\alpha}{\to}\MM_2\overset{\beta}{\to}\MM_3
\end{equation}
by setting $\alpha\colon \mu \mapsto \Psi$ with $\Psi^\xi=\mu(\xi_M)$
and $\beta$ to be the homomorphism
$$
\PP \otimes \g^* \to \bigoplus_I \PP_I \otimes \g^*_I
$$
associated with the restrictions $\PP \to \PP_I$ and $\g^* \to \g_I^*$.
In other words, the $I$th component of $\beta(\Psi)$ is the
$\g_I^*$-component of $\Psi$ restricted to $(\Cx)^I$. Hence, $\Psi$ 
satisfies the second condition of an abstract moment map 
if and only if $\beta(\Psi)=0$, and $\Psi$ is associated with a one-form
if and only if it is in the image of $\alpha$.
Proposition \ref{polynomial} is equivalent to
the sequence \eqref{eq:(2)} being exact.

Denote by $\OO$ and $\EE$, respectively, the algebras of germs of
analytic, resp.\ smooth, functions on $\C^d$ at the origin.
Let $\MMs_i$, where $i=1,\ 2,\ 3$, be the modules
defined similar to $\MM_i$ but in the category of smooth germs at the
origin. Note that $\MMs_i$ are modules over $\EE$.

As before, we have a sequence of homomorphisms
\begin{equation} \labell{exact-smooth}
\MMs_1\overset{\alpha}{\to}\MMs_2\overset{\beta}{\to}\MMs_3.
\end{equation}
To prove the theorem in the smooth category,
it suffices to show that this sequence is exact.

Note that the inclusions $\PP\to\OO$ and $\PP\to\EE$ make
$\OO$ and $\EE$ into $\PP$-modules.
The following lemma is obvious:

\begin{Lemma}
\labell{lemma:tensor}
$\MMs_i=\MM_i\otimes_{\PP}\EE$.
\end{Lemma}

To finish the proof of Theorem \ref{twoform-linear}, we need to recall some 
facts from commutative algebra. Let $B$ be a commutative ring and $A$ a 
sub-ring of $B$. The ring $B$ is said to be {\em flat\/} over $A$ if for 
every exact sequence of $A$-modules
$$
\MM_1\to\MM_2\to\MM_3
$$
the sequence
$$
\MM_1\otimes_A B\to\MM_2\otimes_A B\to\MM_3\otimes_A B
$$
is also exact.

It is known that $\OO$ is flat over $\PP$ (see \cite{Mal},
page 45, Example 4.11) and $\EE$ is flat over $\OO$ (see \cite{Mal},
page 88, Corollary 1.2). This in turn implies that $\EE$ is flat
over $\PP$.
Therefore, the exactness of \eqref{eq:(2)} implies the exactness
of \eqref{exact-smooth}.
\end{proof}

\section{Assignment cohomology}
\labell{sec:cohomology}

In this section we show that the space of assignments $\AA(M)$
on a $G$-manifold
fits as the zeroth space in a sequence of vector spaces
$\HA^*(M)$, called the assignment cohomology.

\subsection{Construction of assignment cohomology}

Let $M$ be a manifold with an action of a torus $G$.
Denote by $P_M$ the set of its infinitesimal orbit type strata
(see Section \ref{sec:def-examples}).
For each stratum $X$, denote by $\g_X$ the infinitesimal stabilizer
of the points of $X$, and let 
$$
V(X) = \g_X^*
$$ 
be the dual space.
Recall that $P_M$ is a partially ordered set, a \emph{poset} for brevity,
with $X \preceq Y$ if  $X$ is contained in the closure of $Y$. 
Denote by $\pi_Y^X \colon V(X) \to V(Y)$ 
the natural projection dual to the inclusion map $\g_Y \subseteq \g_X$
when $ X \preceq Y$.
Recall that the space of assignments is
$$ 
\AA(M)=\{ v \in \prod \limits_{X \in P_M} V(X) \ | \ \pi^X_Y v_X = v_Y 
\text{ for all } X \preceq Y \}.
$$
Every abstract moment map induces an assignment. We will call elements
of $\AA(M)$ \emph{moment assignments}, to distinguish them from 
assignments with other coefficients $V(X)$, which we introduce later.

Define the \emph{assignment cohomology} $\HA^*(M)$ to be the cohomology
of the following cochain complex  which we will denote by $C^*(M;V)$.
A $k$-cochain is a function $\varphi$ that associates to each ordered 
$(k+1)$-tuple $X_0 \preceq \ldots \preceq X_k$ of elements of $P_M$
an element in $V(X_k)$. The differential $d$ is defined by the formula
\begin{eqnarray}
\labell{eq:diff}
d\varphi(X_0,\ldots,X_{k+1})
&=&\sum _{\ell=0}^{k} (-1)^\ell 
  \varphi(X_0,\ldots, \widehat{X_\ell},\ldots, X_{k+1}) \nonumber\\
&&\quad +(-1)^{k+1}\pi^{X_k}_{X_{k+1}}\varphi(X_0,\ldots, X_{k})
,\end{eqnarray}
where, as usual, the hat over $X_\ell$ means that $X_\ell$ is omitted.

\begin{Example}
The zeroth assignment cohomology is simply the space of assignments:
$$ \HA^0(M) = \AA(M).$$
\end{Example}

\begin{Remark}[Functoriality] 
Assignment cohomology is functorial with respect to 
equivariant maps of manifolds:

Let $M$ and $N$ be $G$-manifolds and let $f\colon M \to N$ be a $G$-equivariant
map.   Such a map might not send a stratum in $M$ to a stratum in $N$.
(For example, the function $f(z,w) = z$ from $\C^2$ with the diagonal
circle action to $\C$ with the standard circle action sends the open dense
stratum $\C^2 \ssminus 0$ to the union of strata 
$\{ 0 \} \sqcup \C^\times = \C$.)
However, it does induce a monotone mapping of posets, 
$\tilde{f} \colon P_M \to P_N$, in the following way.
For each stratum $X$ in $M$ there exists 
a unique stratum $Y$ in $N$ with $\g_{X} \subseteq \g_{Y}$ whose closure
contains $f(X)$. (To see this, consider the infinitesimal $\g_X$-action
in $N$. Since $X$ is connected, $f(X)$ is contained in a unique
component of the $\g_X$-fixed point set of $N$. This component is
a $G$-invariant submanifold of $N$. The stratum $Y$ is the open
dense stratum in this component.)
We set $\tilde{f}(X) = Y$. Note that $\g_X \subseteq \g_{\tilde{f}(X)}$
for all $X$. By definition, the pullback map on the cochain complexes
sends a cochain $\varphi \in C^k(N;V)$ to the cochain
$f^*\varphi \in C^k(M;V)$ given by
$$
	(f^*\varphi)(X_0,\ldots,X_k) 
= \pi^{\tilde{f}(X_k)}_{X_k}\varphi(\tilde{f}(X_0),\ldots,\tilde{f}(X_k)).
$$
The map $f^*$ commutes with $d$ and thus induces a pullback map in cohomology.
\end{Remark}

\begin{Theorem} \labell{vanishes}
$\HA^k(M) = 0$ when $k \geq \dim M$ or $k \geq \dim G$.
\end{Theorem}

The proof of Theorem \ref{vanishes} will be an easy application of the 
following alternative 
construction of assignment cohomology. Let $C^k_0(M;V)$ be the space of 
functions $\varphi$ that associate to each ordered $(k+1)$-tuple 
$ X_0 \prec X_1 \prec \ldots \prec X_k$ of \emph{distinct} elements of $P_M$
an element of $V(X_k)$. This can be identified with the subspace of 
$C^k(M;V)$ consisting of those cochains that vanish on $(X_0,\ldots,X_k)$ 
whenever $X_i=X_{i+1}$ for some $i$.  

\begin{Theorem} \labell{C0} 
$C^*_0(M;V)$ is a subcomplex of $C^*(M;V)$, and the inclusion map of complexes
induces an isomorphism in cohomology.
\end{Theorem} 

This result is standard. However, for the sake of completeness, 
we prove it below.

The complex $C^*_0(M;V)$ is much smaller than $C^*(M;V)$ and is more 
convenient to use for explicit calculations. Its disadvantage is that
this complex is not functorial with respect to mappings of posets:
a map $f \colon M \to N$ that send strata to strata sends a tuple 
$X_0 \prec \ldots \prec X_k$ to a tuple 
$f(X_0) \preceq \ldots \preceq f(X_k)$, but the $f(X_j)$ might not be distinct
even if the $X_j$ are. 

\begin{proof}[Proof of Theorem \ref{vanishes}]
By Theorem \ref{C0} it suffices to prove Theorem \ref{vanishes} for the
cohomology of the complex $C^*_0(M;V)$.

The first part of the theorem follows from the fact that
if $X \prec Y$ then $\dim X < \dim Y$.
Therefore, the longest possible tuple
$X_0 \prec \ldots \prec X_k$ of distinct strata
has $k = \dim M$. Thus $C^k_0(M;V) = 0$ for $k > \dim M$.
When $k = \dim M$, the maximal stratum $X_k$ is the open dense
stratum in $M$, on which $V(X_k) = 0$. As a result, $C^{\dim M}_0 (M;V) = 0$.

The second part of the theorem follows from the fact that
if $X \prec Y$ then $\dim \g_X > \dim \g_Y$.
The same argument as before shows that
$C^k_0(M;V) = 0$ whenever $k \geq \dim G$.
\end{proof}

Notice that the proof of the second part of the theorem breaks down 
if the infinitesimal orbit type stratification 
is replaced by the orbit type stratification.  

\begin{proof}[Proof of Theorem \ref{C0}]
Denote by $(X_0^{k_0}, \ldots, X_l^{k_l})$ the tuple
$$
(X_0,\ldots,X_0,\ldots,X_l,\ldots,X_l)
$$
in which each $X_j$ occurs $k_j$ times and the strata $X_j$ are ordered 
and distinct:
$X_0 \prec X_1 \prec \ldots \prec X_l$.
Denote by $n(X_0^{k_0}, \ldots, X_l^{k_l})$ the number of $j$'s such that 
$k_j > 1$; call this number the \emph{fatness} of the tuple.
As is easy to check, the fatness of a $(k+1)$-tuple is no greater than $k$.

For each integer $n \geq 0$ let $C^k_n(M;V)$ be the space of $(k+1)$-cochains 
that are supported on tuples of fatness $n$.  
(This is consistent with the previous definition of $C^k_0(M;V)$.) Then
$$   C^k(M;V) = C^k_0(M;V) \oplus \ldots \oplus C^k_k(M;V)$$
as vector spaces. Set 
$$
C^k_{>0}(M;V) = \bigoplus\limits_{n=1}^k C^k_n(M;V)
.$$
For every nonzero cochain $\varphi$ in this space, there exists a unique 
integer $n$ between $1$ and $k$ and a unique decomposition 
\begin{equation} \labell{varphi1k}
 \varphi = \varphi_n + \varphi_{n+1} + \ldots +\varphi_k
\end{equation}
such that $\varphi_j \in C^*_j(M;V)$ for all $j$ 
and such that $\varphi_n \neq 0$.

An easy computation shows that 
$d(C^k_n (M;V)) \subseteq C^{k+1}_n (M;V) + C^{k+1}_{n+1} (M;V)$
for all $n \geq 1$ and that
$d(C^k_0 (M;V)) \subseteq C^{k+1}_0 (M;V)$.
In particular, $C^*_0(M;V)$ and $C^*_{>0}(M;V)$ are subcomplexes,
and the assignment cohomology splits:
$$ \HA^*(M;V) = \HA^*_0(M;V) \oplus \HA^*_{>0}(M;V).$$
It remains to show that $\HA^*_{>0}(M;V)$ vanishes.

Define a linear map $L \colon C^k(M;V) \to C^{k-1}(M;V)$ by
$$ (L \varphi) (X_0^{k_0} , \ldots , X_l^{k_l}) 
   = \sum_{j=0}^l (-1)^{k_0 + \ldots + k_{j-1}} 
               \varphi(X_0^{k_0}, \ldots , X_j^{k_j+1}, \ldots ,X_l^{k_l}).$$
An explicit computation shows that
\begin{equation} \labell{dL}
 dL\varphi + Ld\varphi = \pi \varphi,
\end{equation}
where 
$$(\pi \varphi) (X_0^{k_0} ,\ldots , X_l^{k_l}) 
       = n(k_0,\ldots,k_l) \varphi(X_0^{k_0},\ldots,X_l^{k_l}).$$
Therefore, $\pi\colon C^*(M;V)\to C^*(M;V)$ is chain homotopic
to zero and, as a consequence, induces the zero map $\pi_*$ in homology.
Denote by $j\colon C^*_{>0}(M;V)\to C^*(M;V)$  the natural inclusion.
It is clear that $\pi j\colon C^*_{>0}(M;V)\to C^*_{>0}(M;V)$ is
an isomorphism. Thus $\pi_* j_*\colon \HA^*_{>0}(M;V)\to \HA^*_{>0}(M;V)$
is an isomorphism. Since $\pi_*=0$, this is possible only when
$\HA^*_{>0}(M;V)=0$. 

\end{proof} 

\subsection{Assignments with other coefficients}
\labell{subsec:other coefficients}
The definitions of assignments and assignment cohomology extend
word--for--word to other systems of coefficients.
A system of coefficients $V$ on the poset $P_M$
is a function that associates a vector space $V(X)$ to each stratum 
$X$ and a linear map $ \pi^X_Y \colon
V(X) \to V(Y)$ to each pair $X\preceq Y$, so that $\pi_Z^Y \pi_Y^X =
\pi_Z^X$ whenever $X \preceq Y \preceq Z$ and $\pi_X^X = \id$.  We define
a differential complex $(C^*(M;V),d)$ as before and denote its cohomology
by $\HA^*(M;V)$.

A morphism $V_1 \to V_2$ of two systems of coefficients consists of
a linear map $V_1(X) \to V_2(X)$ for each $X \in P_M$, such that
the squares 
$$ 
	\begin{array}{ccc}
	V_1(X) & \longrightarrow & V_2(X) \\
	\downarrow & & \downarrow \\
	V_1(Y) & \longrightarrow & V_2(Y) 
	\end{array}
$$
commute for all $X \preceq Y$. 
The systems of coefficients $V$ on $P_M$ form a category. The assignment 
cohomology groups, $\HA^n(M;V)$, and in particular the space of 
assignments $\AA(M;V)$, are functorial in $V$.

\begin{Remark} \labell{derived functors}
One can think of the poset $P_M$ as a category in which there exists
a single arrow $Y \to X$ whenever $X \preceq Y$. A system of coefficients
is a contra-variant functor from $P_M$ to the category of vector spaces,
and a morphism of systems of coefficients is a natural transformation.
The space of assignments is the inverse limit:
\begin{equation} \labell{liminv}
 \AA(M;V) = {\liminv}_{X \in P_M} V(X),
\end{equation}
which is a functor in $V$.\footnote{\labell{foot:order} 
Note that the convention on the direction of morphisms commonly used
to turn a poset into a category is opposite of the one employed in this
paper. Alternatively, the poset of strata is sometimes given
the order inverse of the one used above. However, the only essential 
point in the choice of directions of morphisms is that a poset should be
made into a category so that $V$ becomes a \emph{contra-variant} functor.}

A system of coefficients $V$ can be viewed as a \emph{pre-sheaf} on 
the category $P_M$,
in the sense of \cite[Expose I, Definition 1.2]{SGA4}. 
(Also see \cite[Ch.~I \S 2 and \S 4]{Mo}.)
The assignments are the global sections:
\begin{equation} \labell{Gamma}
 \AA(M;V) = \Gamma(P_M;V) ,
\end{equation}
and the assignment cohomology is equal to the cohomology of this presheaf. 
More explicitly, the cohomology groups of the presheaf are defined to be 
the derived functors of the global section functor \eqref{Gamma} (as a functor 
in $V$); equivalently, the cohomology groups are the derived functors 
of the inverse limit functor \eqref{liminv}.
The fact that these derived functors are the same as the cohomology of
the complex $C^*(M;V)$ is Proposition 6.1 in \cite[ch.~II]{Mo}.
\end{Remark}

For some applications it is beneficial to work with systems of 
coefficients with values in categories other than the category
of vector spaces. Let us illustrate this by some examples.

\begin{Example} \labell{isotropy AA}  
Assignments with values in the functor
\begin{equation} \labell{virtual reps}
 V(X) = \{ \text{ equivalence classes of representations of $\g_X$ } \}.
\end{equation}
are called \emph{isotropy assignments}. A $G$-action gives rise
to a canonical isotropy assignment defined as follows:
to each $X$ associate the isotropy representation of $\g_X$ on $T_pM$,
$p \in X$. In a similar manner,
any $G$-equivariant vector bundle over $M$
gives rise to an isotropy assignment.
\end{Example}

\begin{Example}
For a symplectic manifold with a Hamiltonian torus action,
its X-ray is, roughly speaking, the direct sum of the isotropy assignment 
and the moment map assignment. The notion of X-rays, introduced 
in \cite{tolman}, is central to the study of Hamiltonian torus actions.
See \cite{tolman,hamiltonian,defone,metzler1,metzler2}.
\end{Example}

A broad class of systems of coefficients can be obtained by the following 
construction.
Consider the category of sub-algebras of $\g$ with morphisms given
by inclusion maps. Let $V'$ be any functor from this category to the
category of vector spaces (or modules, abelian groups, etc.).  Then we get
a system of coefficients with $V(X) = V'(\g_X)$.
For instance, we get moment assignments from $V'(\h) = \h^*$,
and we get isotropy assignments from 
$V'(\h) = \{ \text{ virtual representations of } \h \}$. 

\begin{Example}
\labell{exam:min2}
Assume that $M$ has a unique minimal stratum, $X_0$.
This is the case, for instance, when $M$ is a vector space
on which $G$ acts linearly. Then $\HA^0(M;V)=V(X_0)$ and
$\HA^k(M,V)=0$ for all $k>0$. Indeed, the operator
$$
	(Q\varphi) (Y_1,\ldots,Y_k) = \varphi(X_0,Y_1,\ldots,Y_k) 
$$
satisfies $dQ\varphi+Qd\varphi=\varphi$, hence it is a homotopy
operator for the complex $C^*(M;V)$. (Also see Remark \ref{rmk:min}.)
\end{Example}

Theorem \ref{C0} holds for assignment cohomology with many other
systems of coefficients. For example,
the theorem clearly holds whenever $V$ takes values in the category of
vector spaces. When Theorem \ref{C0} applies,
we also have the following variant of Theorem \ref{vanishes}:

\begin{Proposition}
\begin{enumerate}
\item Let $V$ be such that $V(X)=0$ for the open stratum $X$.
Then $\HA^k(M;V)=0$ when $k\geq \dim M$. 
\item Let $V$ be obtained as the pull-back of a functor on the 
sub-algebras of $\g$ which vanishes on the zero sub-algebra. 
Then $\HA^k(M;V)=0$ when $k\geq \dim G$.  
\end{enumerate}
\end{Proposition}
The proof of this fact is entirely similar to the proof of Theorem
\ref{vanishes}.

\subsection{Assignment cohomology for pairs}

Let $N$ be a subset of $M$ which is a union of strata.
Define the \emph{relative} assignment cohomology $\HA^*(M,N;V)$ 
to be the cohomology of the sub-complex $C^*(M,N;V)$
of $C^*(M;V)$ formed by those cochains which vanish on all
$(k+1)$-tuples $X_0 \preceq \ldots \preceq X_k$ in which all of 
the strata $X_j$ are in $N$.

\begin{Theorem} \labell{thm:coh}
There is a long exact sequence
\begin{equation}
\labell{eq:long-exact}
\ldots\to \HA^*(M,N; V)\to \HA^*(M; V)\to \HA^*(N)\to \HA^{*+1}(M,N; V)
\to \ldots ,
\end{equation}
where the connecting homomorphism is given by the standard formula.
\end{Theorem}

\begin{proof}
The theorem is an immediate consequence of the fact
that the sequence of complexes 
\begin{equation} \labell{short exact}
0\to C^*(M,N;V)\to C^*(M;V)\to C^*(N;V)\to 0
\end{equation}
is exact. To prove the exactness of \eqref{short exact},
note that the space $C^*(M,N;V)$ is the kernel of the restriction map
$C^*(M;V) \to C^*(N;V)$ by its definition.
To see that the restriction map is onto, note that every cochain
$\varphi$ in $C^*(N;V)$ can be extended to a cochain 
$\tilde{\varphi}$ in $C^*(M;V)$ by declaring 
$\tilde{\varphi}(X_0,\ldots,X_k)$ to be zero whenever not all of $X_j$'s
are in $N$.
\end{proof}

An alternative proof of Theorem \ref{thm:coh} in the case where $N$
is open is given in Remark \ref{rmk:pf-long-exact}.

\begin{Remark}
\labell{rmk:exact-seq'} 
It is not clear if there is a way to define relative assignment 
cohomology so that the sequence \eqref{eq:long-exact} is exact 
in the case where $N \subseteq M$ is $G$-invariant but not 
a union of strata.
 For instance, consider a $G$-manifold $B$ and let $G$ 
act on $M = B \times [0,1]$ by acting on the first factor.
Let $N = B \times \{ 0, 1 \}$. Since every stratum in $M$
meets $N$, it seems reasonable to set $\HA^*(M,N)=0$.
Then an exact sequence \eqref{eq:long-exact} 
would give an isomorphism between $\AA(M)$
and $\AA(N)$. However, these spaces are not isomorphic; in fact,
$\AA(N)$ is isomorphic to $\AA(M) \oplus \AA(M)$.
\end{Remark}

As is with other ``cohomology theories'', an exact sequence of
coefficients gives rise to a long exact sequence in assignment cohomology:

\begin{Theorem} \labell{short-long}
A short exact sequence of systems of coefficients,
\begin{equation} \labell{V123}
0 \to V_1 \to V_2 \to V_3 \to 0,
\end{equation}
induces a long exact sequence in cohomology,
\begin{equation} \labell{HkV123}
\to \HA^k(M,N;V_1) \to \HA^k(M,N;V_2) \to \HA^k(M,N;V_3) 
\stackrel{\delta}{\to} \HA^{k+1}(M,N;V_1) \to 
\end{equation}
\end{Theorem}

\begin{proof}
The short exact sequence \eqref{V123} naturally induces a short exact
sequence of complexes,
$$ 0 \to C^*(M,N;V_1) \to C^*(M,N;V_2) \to C^*(M,N;V_3) \to 0,$$
and hence the long exact sequence \eqref{HkV123} in cohomology.
\end{proof}

\begin{Remark} \labell{M,N}
Suppose that $N \subseteq M$ is an \emph{open} subset which is a union 
of strata. Note that being open is equivalent to the following condition:
\begin{equation} \labell{N complete}
\text{ for every pair } X \prec Y \text{ of strata in } M,
\text{ if } X \subseteq N \text{ then also } Y \subseteq N.
\end{equation}
Then the relative assignment cohomology is equal to the ordinary  
assignment cohomology with a different system of coefficients:
$$ 
	\HA^n(M,N;V) = \HA^n(M;V_{M/N})
$$
where $V_{M/N}$ is the system of coefficients given by
$$ 
	V_{M/N}(X) = 
	\begin{cases}
		V(X) & \text{if $X \not \subseteq N$},\\
		0 & \text{otherwise}
	\end{cases} 
$$
for any stratum $X$ in $M$, with the projection maps
$$
(\pi_{M/N})^X_Y =
\begin{cases}
\pi^X_Y & \text{if both $X$ and $Y$ are not in $N$,}\\
0 & \text{otherwise.}
\end{cases}
$$
The compatibility condition,
\begin{equation} \label{compatibility}
(\pi_{M/N})^Y_Z \circ (\pi_{M/N})^X_Y = (\pi_{M/N})^X_Z \text{ whenever } 
X \preceq Y \preceq Z, 
\end{equation}
which is required for $V_{M/N}$ to form a system of coefficients, 
is satisfied if and only if $N$ meets the requirement \eqref{N complete}.
\end{Remark}

\begin{Remark}
\labell{rmk:pf-long-exact}
If $N \subseteq M$ is a union of strata and is open, the long exact sequence 
for the pair $(M,N)$ in Theorem \ref{thm:coh} follows from the long 
exact sequence for coefficients in Theorem \ref{short-long}.  To see this, 
recall that $\HA^n(M,N;V) = \HA^n(M;V_{M/N})$ as explained in 
Remark \ref{M,N}.   Furthermore, we set $V_N(X)=V(X)/V_{M/N}(X)$
with $(\pi_N)^X_Y=\pi^X_Y$ if $X$ and $Y$ are both in $N$ and
$(\pi_N)^X_Y=0$ otherwise.  Then we have $\HA^*(M;V_N)=\HA^*(N;i^*V)$, 
where $i \colon N \to M$ is the inclusion map.
The sequence of systems of coefficients
$ 0 \to V_{M/N} \to V \to V_N \to 0 $
is exact. By Theorem \ref{short-long}, this sequence gives rise
to the long exact sequence which coincides with the sequence
\eqref{eq:long-exact} of Theorem \ref{thm:coh}.
\end{Remark}

\begin{Remark} \labell{delta}
Relative assignment cohomology is a sequence of functors 
$$
	V \mapsto \HA^n(M,N;V)
$$ 
from the category $\calC_M$ of systems of coefficients 
on $P_M$ to the category of vector spaces.
This sequence, together with the maps $\delta$ of \eqref{HkV123},
form a \emph{$\delta$-functor}. (Essentially, this means
that short exact sequences in $\calC_M$ induce long exact sequences
in cohomology, as in Theorem \ref{short-long}.  See \cite{lang}.)

In the non-relative case $N=\emptyset$ this $\delta$-functor is 
\emph{universal}, i.e, the functors $V \mapsto \HA^n(M;V)$ 
are the derived functors of the assignment functor 
$V \mapsto \AA(M;V)$; see Remark \ref{derived functors}.
This remains true in the relative case if $N \subseteq M$ is a union of 
strata and is open: the relative assignment cohomology functors 
$V \mapsto \HA^n(M,N;V)$ are then the derived functors of the relative 
assignment functor $V \mapsto \AA(M,N;V)$ which associates to each $V$
the space of assignments that vanish on $N$. However, for a general $N$,
it is not clear if these functors are universal or,
equivalently, whether or not $V\mapsto \HA^*(M,N;V)$ are the derived
functors for $V\mapsto \AA(M,N;V)$.
\end{Remark}

\begin{Remark}
\labell{rmk:min}
The poset of strata $P_M$ does not in general satisfy
the following condition which is routinely required in some sources 
(e.g., \cite{Je}, \cite{Mas}, and \cite{Ru}):
\begin{equation} \labell{eq:condP}
	\begin{array}{c}
		\text{for any } X \in P_M \text{ and } Y \in P_M 
		\text{ there exists } Z \in P_M\\
 		\text{ such that } Z \preceq X \text{ and } Z \preceq Y.
	\end{array}
\end{equation}
(The reader should keep in mind that our order convention is opposite
of the standard one; see footnote \ref{foot:order}.)
This condition is met, for example, when $P_M$ has a minimal element,
i.e., an $X_0 \in P_M$ such that $X_0 \preceq X$ for all $X \in P_M$.
(Equivalently $X_0$ is a stratum which is contained in
the closure of every stratum $X$.)
With the condition \eqref{eq:condP}, the poset $P$ is called
a \emph{directed set}. Under this condition, 
an inverse system $V$ of finite--dimensional vector spaces is 
automatically flabby (see, e.g., \cite{Je} and \cite{Ru}) and 
$\liminvk V=0$ for all $k>0$.
This generalizes Example \ref{exam:min2}.
\end{Remark}

\subsection{Examples of calculations of assignment cohomology}

The following simple example shows that relative assignment
cohomology  can be non-trivial in degrees greater than zero.

\begin{Example}
\labell{exam:relative}
Let $M=\CP^2$ and $G$ be the torus $\T^2$ acting on $M$ as in Example 
\ref{CP2-T2}, and $N=M^G$. Then $\AA(M,N)=0$, $\AA(N)=(\g^*)^3$ is 
six-dimensional, and $\dim\AA(M)=3$. Furthermore,
$\HA^{*>0}(M)=\HA^{*>0}(N)=0$. Thus \eqref{eq:long-exact} turns into
the exact sequence
$$
	0 \to \AA(M) \to \AA(N)\to \HA^1(M,N) \to 0 ,
$$
where $\dim \HA^1(M,N)=3$, as can also be checked by a direct calculation.
\end{Example}

\begin{Example}[Assignment cohomology for toric varieties]
Let $M$ be a compact smooth K\"ahler toric manifold of complex dimension $n$
with moment map $\Psi \colon M \to \g^*$.
(See Examples \ref{toric:1} and \ref{toric:2}.)
Recall that the poset $P_M$ of orbit type strata is isomorphic to the poset 
of faces $\Psi(X)$ of a simple polytope $\Psi(M)$, and for each stratum $X$,
$$
	\dim_\C X = \dim \Psi(X) = n - \dim \g_X^*.
$$
We will work with the system of coefficients $V(X) = \g_X^*$.
The zeroth assignment cohomology is the space of assignment
which was computed in Example \ref{toric:1}. Namely,
$$
\HA^0(M;V)=\bigoplus_Y\g_Y^*,
$$
where the summation is over all strata $Y$ with $\dim \g_Y=1$. These
strata correspond to the $(n-1)$-dimensional faces of $\Psi(M)$. In
particular,
$$
	\dim \HA^0(M;V) = \text{the number of facets of $\Psi(X)$}.
$$

Let us prove that the higher cohomology groups vanish:
\begin{equation} \labell{higher vanish} 
	\HA^k(M;V) = 0 \quad \text{for all $k \geq 1$.}
\end{equation}
By Theorem \ref{C0} it is enough to work with the complex $C^*_0(M;V)$.
For a closed cochain $\varphi \in C^k_0(M;V)$, $k \geq 1$, we will find
a primitive $(k-1)$-cochain $\psi \in C^{k-1}_0(M;V)$, i.e., 
a cochain $\psi$ such that $d\psi=\varphi$.

Let $X_0 \prec \ldots \prec X_{k-1}$ be any ordered $k$-tuple of distinct
strata. Recall that the natural map
\begin{equation} \labell{projections}
	V(X_{k-1}) \stackrel{\oplus \pi^{X_{k-1}}_{X_k}}{\longrightarrow}
	   \bigoplus\limits_{X_k}V(X_k),
\end{equation}
where $X_k$ is such that $\codim X_k=1$ and $X_{k-1}\preceq X_k$,
is a linear isomorphism. Therefore, to define the value 
$\psi(X_0,\ldots,X_{k-1})$, which is an element of $V(X_{k-1})$, 
it is enough to specify the projections of these elements to all of the
spaces $V(X_k)$ with $X_k$ as above.
We require these projections to be
\begin{equation} \labell{pi psi}
	\pi^{X_{k-1}}_{X_k} \psi(X_0,\ldots,X_{k-1}) 
	= (-1)^k\varphi(X_0, \ldots, X_k).
\end{equation}
Let us show that $d\psi=\varphi$.
Set $\varphi' = \varphi - d\psi$. Note that
$\psi(X_0,\ldots,X_{k-1})=0$ when $\codim X_{k-1}=1$. Then it follows
from the definition \eqref{pi psi} of $\psi$
and the definition \eqref{eq:diff} of the differential that $\varphi'$ 
vanishes on all tuples $X_0 \prec \ldots \prec X_k$ in which 
$\codim X_k = 1$. Again, by \eqref{eq:diff},
$d\varphi'(X_0,\ldots,X_{k+1}) = \pi^{X_k}_{X_{k+1}} \varphi'(X_0,\ldots,X_k)$
for all tuples $X_0 \prec \ldots \prec X_k \prec X_{k+1}$ in which $\codim X_{k+1} = 1$.
Since $d\varphi'=0$, and since \eqref{projections} is a linear isomorphism,
this implies that $\varphi' \equiv 0$.  Hence, $\varphi = d\psi$ is exact.
\end{Example}

We now give an example of a manifold which has a non-trivial 
(absolute, not relative) first assignment cohomology. 

\begin{Example}
Let $M = S^2 \times S^2 \times S^2$, and let $G = S^1 \times S^1$
act by 
$$
	(a,b) \cdot (u,v,w) = (a \cdot u, b \cdot v, ab\inv \cdot w)
$$
where on the right the dot denotes the standard $S^1$ action 
on $S^2$ by rotations. The moment assignments can be drawn as pictures
showing the moment map images of the orbit type strata (the ``x-ray").
Such a picture is shown in Figure \ref{fig:example}. 
Notice that this picture is two--, not three--, dimensional.
This arrangement
can be moved around as long as the edges are shifted but not rotated.
An assignment is therefore determined by the location of the bottom left
vertex and the lengths of the three edges coming out of it. Therefore, 
\begin{equation} \labell{five}
	\dim \HA^0(M;V) = \dim \AA(M;V) = 5.
\end{equation}
We will find the dimension of the first assignment cohomology
space by using the Euler characteristic of the complex $C^*_0(M;V)$ 
of Theorem \ref{C0}.  We have
\begin{equation} \labell{euler}
	\begin{array}{ccl}
		\dim \HA^0(M;V) - \dim \HA^1(M;V) 
			& = & \sum_k (-1)^k \dim C^k_0(M;V) \\
		& = & \dim C^0_0(M;V) - \dim C^1_0(M;V) 
\end{array}
\end{equation}
because $C^k_0(M;V) = 0$ for all $k \geq 2$ (see also Theorem \ref{vanishes}).
A $0$-cochain associates to each vertex an element of a two--dimensional
space and to each edge an element of a one--dimensional space. Therefore,
\begin{equation} \labell{28}
	\begin{array}{ccl}
		\dim C^0_0 (M;V) & = & 2 (\text{number of vertices}) 
				+ (\text{number of edges}) \\
 			& = & 2 \cdot 8 + 12 \\
 			& = & 28.
	\end{array}
\end{equation}
A $1$-cochain associates an element of a one dimensional space
to each pair consisting of a vertex and an edge coming out of it.  Therefore
\begin{equation} \labell{24}
	\dim C^1_0 (M;V) = 24.
\end{equation}
Substituting \eqref{28}, \eqref{24}, and \eqref{five} in \eqref{euler},
we get
$$
	5 - \dim \HA^1(M;V) = 28 - 24,
$$
hence
$$
	\dim \HA^1(M;V) = 1.
$$
\end{Example}

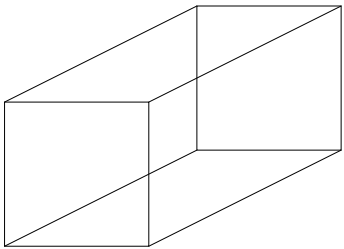
\begin{figure}
\begin{center}
\setlength{\unitlength}{0.00083333in}
\begingroup\makeatletter\ifx\SetFigFont\undefined%
\gdef\SetFigFont#1#2#3#4#5{%
  \reset@font\fontsize{#1}{#2pt}%
  \fontfamily{#3}\fontseries{#4}\fontshape{#5}%
  \selectfont}%
\fi\endgroup%
{\renewcommand{\dashlinestretch}{30}
\begin{picture}(2124,1539)(0,-10)
\path (12,12) (912,12) (912,912) (12,912) (12,12)
\path (1212,612) (2112,612) (2112,1512) (1212,1512) (1212,612) 
\path (12,12) (1212,612)
\path (912,12) (2112,612)
\path (912,912) (2112,1512)
\path (12,912)(1212,1512) 
\end{picture}
}
\end{center}
\caption{A manifold with nonzero first assignment cohomology}
\labell{fig:example}
\end{figure}

\subsection*{Generalizations of assignment cohomology}
The results and definitions of this section can be generalized or altered 
in many natural ways. For instance, the assignment cohomology can be
defined for an arbitrary system of coefficients with values in an
abelian category on an arbitrary poset. In particular, in a more 
geometrical realm, the infinitesimal orbit type 
stratification can be replaced by the orbit type stratification and $V$ 
can then be the pull-back of a contra-variant functor on subgroups of 
$G$. Furthermore, instead of working with functors with values
in finite--dimensional vector spaces one may consider functors
with values in abelian groups or graded vector spaces or rings.
In fact, such functors do arise arise in the study of symplectic 
manifolds, and $\AA(M;V)$ can be viewed as a repository containing many 
of the invariants of the action.  One important example is the isotropy 
assignments of Example \ref{isotropy AA} above.

Moreover, most of the results of this section extend with obvious 
modifications to actions of finite or compact \emph{non-abelian} groups.
For example, we can take $V$ to be the 
pull-back of a functor on the sub-algebras of $\g$ which is invariant 
under conjugations. However, the space of moment map assignments 
for actions of non-abelian groups (Remark \ref{rmk:non-abelian})
does not arise as the zeroth cohomology groups of this type.
In the non-abelian case, the space of moment map assignments 
does not seem to be associated with a functor on the poset of strata.
Instead, to obtain this space one should work with the singular 
foliation of $M$ given by the decomposition of $M$ according to actual 
stabilizers, but not just their conjugacy classes. This renders a 
correct generalization of assignment cohomology to actions of 
non-abelian groups much less straightforward. 

\begin{Remark}
The assignment cohomology appears to be related to Bredon's  equivariant
cohomology, \cite{Bredon}, and, perhaps, to Borel's equivariant 
cohomology with twisted coefficients.  The nature and explicit form of
these relations are, however, unclear to the authors.
\end{Remark}

\end{document}